\documentclass[]{article}

\usepackage{amssymb}
\usepackage{amsmath}
\usepackage{graphicx}
\usepackage{epsfig}
\usepackage{float}
\usepackage{subfigure}
\usepackage{psfrag}
\usepackage[utf8x]{inputenc}
\usepackage{color}
\usepackage{array}
\usepackage{dsfont}
\usepackage{algorithm,algorithmic}
\usepackage{tikz}
\usepackage{listings}
\usepackage{epsfig}
\usetikzlibrary{arrows,decorations.pathmorphing,backgrounds,fit,positioning,shapes.symbols,chains,shapes}
\usepackage{multirow}
\usepackage{booktabs}
\usepackage{a4wide}

\definecolor{bleuONERA}{RGB}{16,97,169}
\definecolor{grisONERA}{RGB}{64,64,66}

\newtheorem{remark}{Remark}

\providecommand{\mimo}[0]{\textbf{MIMO}~}

\providecommand{\miso}[0]{\textbf{MISO}~}

\providecommand{\lti}[0]{\textbf{LTI}~}

\providecommand{\ie}[0]{\emph{i.e.}~}

\providecommand{\eg}[0]{\emph{e.g.}~}

\DeclareMathOperator*{\rank}{\mathbf{rank}}

\providecommand{\abs}[1]{\left\lvert #1 \right\rvert} %
\providecommand{\norml}[1]{\left\lVert #1 \right\rVert} %
\providecommand{\norm}[1]{|| #1 ||} %
\providecommand{\x}[0]{\mathbf{x}} %
\providecommand{\xr}[0]{\mathbf{\hat{x}}} %
\renewcommand{\u}{\mathbf{u}} %
\providecommand{\y}[0]{\mathbf{y}} %
\providecommand{\z}[0]{\mathbf{z}} %
\providecommand{\zr}[0]{\mathbf{\hat{z}}} %
\providecommand{\lv}[0]{\mathbf{l}} %
\providecommand{\rv}[0]{\mathbf{r}} %
\providecommand{\wv}[0]{\mathbf{w}} %
\providecommand{\vv}[0]{\mathbf{v}} %

\providecommand{\Hrealr}[0]{\mathcal{\hat{S}}} %
\providecommand{\Ar}[0]{{\hat{A}}} %
\providecommand{\Br}[0]{{\hat{B}}} %
\providecommand{\Cr}[0]{{\hat{C}}} %
\providecommand{\Dr}[0]{{\hat{D}}} %
\providecommand{\Nr}[0]{{\hat{N}}} %
\providecommand{\Fr}[0]{{\hat{F}}} %

\providecommand{\Hreal}[0]{\mathcal{S}} %
\providecommand{\Htran}[0]{\mathbf{H}} %
\providecommand{\E}[0]{{E}} %
\providecommand{\A}[0]{{A}} %
\providecommand{\B}[0]{{B}} %
\providecommand{\C}[0]{{C}} %
\providecommand{\D}[0]{{D}} %

\providecommand{\LL}[0]{{\mathds L}} %
\providecommand{\sLL}[0]{{\mathds L_\sigma}} %

\providecommand{\Htwo}[0]{{\mathcal{H}_{2}}} %
\providecommand{\Hinf}[0]{{\mathcal{H}_{\infty}}} %
\providecommand{\Linf}[0]{{\mathcal{L}_{\infty}}} %

\providecommand{\Cplx}[0]{\mathbb{C}} %
\providecommand{\Real}[0]{\mathbb{R}} %
\providecommand{\matrixtwo}[4]{ \left[\begin{array}{cc} #1 & #2 \\ #3 & #4 \end{array}\right] } %
\providecommand{\vectortwo}[2]{ \left[\begin{array}{c} #1 \\ #2 \end{array}\right] } %
\providecommand{\vectortwoT}[2]{ \left[\begin{array}{cc} #1 & #2 \end{array}\right] } %
\begin{document}


\title{Mixed interpolatory and inference non-intrusive reduced order modeling with application to pollutants dispersion}

\author{C. Poussot-Vassal, T. Sabatier, C. Sarrat and P. Vuillemin
\thanks{Universit\'e de Toulouse, F-31055 Toulouse, France. Contact:  \texttt{charles.poussot-vassal@onera.fr}}
}
\maketitle


\begin{abstract}
    On the basis of input-output time-domain data collected from a  complex simulator, this paper proposes a constructive methodology to infer a reduced-order linear, bilinear or quadratic time invariant dynamical model reproducing the underlying phenomena. The approach is essentially based on linear dynamical systems and approximation theory. More specifically, it sequentially involves the interpolatory Pencil and Loewner framework, known to be both very versatile and scalable to large-scale data sets, and  a linear least square problem involving the raw data and reduced internal variables. With respect to intrusive methods, no prior knowledge on the operator is needed. In addition, compared to the traditional non-intrusive operator inference ones, the proposed approach alleviates the need of measuring the original full-order model internal variables. It is thus applicable to a wider application range than standard intrusive and non-intrusive methods. The rationale is successfully applied on a large eddy simulation of a pollutants dispersion case over an airport area involving multi-scale and multi-physics dynamical phenomena. Despite the simplicity of the resulting low complexity model, the proposed approach shows satisfactory results to predict the pollutants plume pattern while being significantly faster to simulate.
\end{abstract}

\section{Introduction}
\label{sec:intro}
\subsection{Intrusive and non-intrusive reduced order model context}

Dynamical models are central in many engineering fields as they serve for simulation, analysis, optimisation and control. This is even emphasized when considering very complex and critical systems or phenomena for which a deep attention and understanding are needed. These considerations may be further motivated by economical, societal or ecological reasons. In most cases, solutions involving dedicated computer-based software are largely preferred by engineers and researchers to first reduce development costs and time, and second, to improve their understanding of the system under consideration. 

As these systems are usually grounded on accurate complex computationally demanding multi-physics large-scale dynamical models, their approximation by an (accurate) low complexity surrogate dynamical model is then a cornerstone for further advanced developments. This is the purpose of reduced order model (\textbf{ROM}) construction  which can be beneficial in the many-query processes (optimisation, control...). Recent applications can be found in many areas, \eg micro-structures \cite{Vakilzadeh:2020}, fluid mechanics \cite{Willcox:2005,Zimmermann:2014,Star:2021}  or chemistry \cite{Brown:2020}. In fact, many field of studies require physical models with the resolution of Partial Differential Equation (\textbf{PDE}) which have a high computational cost. For example, Computational Fluid Dynamics (\textbf{CFD}) require not only High Performance Computational (\textbf{HPC}) resources but also a lot of preparation work for grid meshing. Concerning, high fidelity atmospheric modeling, based on the resolution of Navier \& Stokes equations, the physics involved is quite complex and the spatio-temporal resolution needs to be very high in order to capture the pollutants dispersion at the local scale (resolution of few meters and one second time step) together with all physical interactions with the environment (buildings, ground, turbulence, radiation, chemistry...). The simulation are usually made using \textbf{HPC} resources and do not enable parametric sensitivity tests for example. To that end, reduced models can be useful in order to perform numerous simulations. The aim of this paper is to show a proof of concept that this type of physical models can be reduced in order to have an approximation fast and easy to run.

An overview of the dynamical model reduction research can be found in books \cite{AntoulasBook:2005,AntoulasBook:2020} or  surveys \cite{AntoulasSurvey:2016,BennerSurvey:2017}. In addition to theoretical works, dedicated (linear) model reduction and approximation numerical tools such as \cite{MORLAB:2017}, \cite{Chebfun} or \cite{mortoolbox} from the authors are also being developed, addressing the closely related computational issues. Reader interested in  theoretical, computational and application insights may also refer to the monograph \cite{PoussotHDR:2019}. Following the definition proposed by \cite{Peherstorfer:2015} and \cite{BennerInfer:2020}, two broad \textbf{ROM} construction can be considered: the \emph{intrusive} and the \emph{non-intrusive} ones. 

Traditional model reduction tools usually refer to \emph{intrusive} methods as they require the original dynamical model mathematical description (\eg operators) to construct the \textbf{ROM}. The most efficient approaches then consist in projecting the operators onto specific subspaces. The definition and choice of the projector and more specifically on its spanned space, is one of the critical topic in this research field, both theoretically and numerically. In the linear case, one can mention the balanced truncation (\textbf{BT}) which projects onto the most observable and controllable subspaces (see \eg \cite{Moore:1981}) or the iterative rational Krylov algorithm (\textbf{IRKA}) for which the spanned subspace is related to the most input-output energetic transfer, minimising some dynamical systems $\Htwo$-norm (see \cite{GugercinSIAM:2008,VanDooren:2008,VanDooren:2010} and a extensions in the parametric case in \cite{BaurSIAM:2011}). Their nonlinear counterpart also exist, as presented in a series of papers addressing bilinear, quadratic or polynomial models (see \eg \cite{Baur:2014} and references therein for more details). 

On the other side, \emph{non-intrusive} methods refer to approaches where only input to state or input to output data are available. Within this category, one can mention input-output frequency-domain data-driven model approximation in the Loewner framework (\textbf{LF}), proposed in its linear version in  \cite{Mayo:2007,AntoulasChapSpringer:2010} and extended to bilinear in  \cite{AntoulasBilinSIAM:2016,GoseaLA:2018} and very recently to time-domain bilinear in \cite{Karachalios:2019}. The \textbf{LF}, as the \textbf{IRKA} (or transfer function \textbf{IRKA}), is an interpolatory method dedicated to data-driven cases where frequency response is accessible, instead of models. One strong property is that the Loewner pencil encodes the minimal rational order of the underlying system. Similarly, identification techniques involving frequency or time-domain data are also proposed through the pencil method in \cite{HoKalman:1966,IonitaPencil:2012}, which is closely related to the Loewner matrices. Generic vector fitting (\textbf{VF}) approaches of \cite{Gustavsen:1999} or adaptive Anderson Antoulas (\textbf{AAA}) proposed in \cite{NakatsukasaAAA:2018} and extended to quadratic outputs in  \cite{GoseaQoutAAA:2020} are also considered as input-output frequency-driven model approximation methods. Dynamic modes decomposition (\textbf{DMD}) proposed in \cite{Proctor:2016} and extended in \cite{GoseaDMD:2020} is an input-state-output time-domain data-driven \textbf{ROM} construction method based on the resolution of a least square problem. Finally, as a \emph{non-intrusive} approach, one also mention the operator inference one proposed in \cite{Peherstorfer:2015} and extended in \cite{BennerInfer:2020}. It allows  inferring a \textbf{ROM} from data collected in the time-domain. A strong property of this last approach is that, as shown in  \cite{Peherstorfer:2015}, it is proved to be equivalent to \emph{intrusive} methods when enough data are collected and under time domain convergence assumptions. Note that this approach has been deployed on a complex combustion process in \cite{McQuarrie:2020}. However, one major limitation of inference and dynamic modes decomposition stands in the need of collecting the internal state time series, which in many cases is impossible because \emph{(i)} it is simply not available or \emph{(ii)} due to storage limitations. 

\subsection{Contributions}

As illustrated later on in the paper, these points are one of the justifications of this work. More specifically, in this paper, we consider the use-case where we have a complex simulator (\eg a Large Eddy Simulation, \textbf{LES}), driven by accessible time-domain input $u(t_k)$ and output $z(t_k)$ data at constantly sampled time instants $\{t_k\}_{k=1}^N$. Importantly, the internal states information is unknown. The objective is, on the basis of these data, to infer a linear, bilinear and/or quadratic dynamical \textbf{ROM} that accurately reproduces the original raw data. The efficiency of this new inferring process is illustrated considering data collected on a \textbf{LES} of pollutants dispersion, showing its efficiency in simulation and applicable to predict the behavior of the plume over a given area in response to pollutant emissions.

\subsection{Paper structure}

The proposed main contribution, being the time-domain input-output mixed interpolatory and operator inference \textbf{ROM} construction process, is detailed in section \ref{sec:inference}. This section typically merges different model approximation and reduction methods recalled in the introduction, for which details are not given to preserve readability. Section \ref{sec:appli} then applies the proposed approach on a very complex set of data obtained from a \textbf{LES} performed with an atmospheric research model to predict the pollutants dispersion over an airport area. The approach shows to very well perform in restitution and opens  perspectives to further developments, sketched in section \ref{sec:conlu}.

\section{Main result: mixed interpolatory and operator inference reduced order model construction from time-domain data}
\label{sec:inference}
The objective of the proposed method is to infer a reduced order dynamical model that accurately reproduces the response of the original full order one, issued from any complex simulator. The specificity in the considered setup is that we consider that we have access to the inputs and outputs (time-domain) raw data only. The overview of the proposed process is first presented in section \ref{ssec:overview}. Section \ref{ssec:data} shows how the data are collected and how preliminary properties of the system can be  extracted, prior any modeling. Then section \ref{ssec:pencil} presents the \emph{pencil method} step, allowing to construct a complex linear dynamical model directly from the time-domain data. Such a model is then referred to as full order model (\textbf{FOM}). A \textbf{ROM} is then computed in section \ref{ssec:rom} by means of an interpolatory method, here in the \emph{Loewner framework}. This latter allows reducing  the dimensionality of the complete model. The loss of accuracy induced by the linear model structure of the \textbf{FOM} construction followed by the reduction step (and discussed later) is then corrected and enriched in section \ref{ssec:infer} by a \emph{model inference} approach exploring different model structures, being either linear or (nonlinear) bilinear or quadratic first-order difference equations, referred to as \textbf{L-ODE}, \textbf{B-ODE} or \textbf{Q-ODE} (or \textbf{QB-ODE}).

\subsection{Procedure at a glance}
\label{ssec:overview}

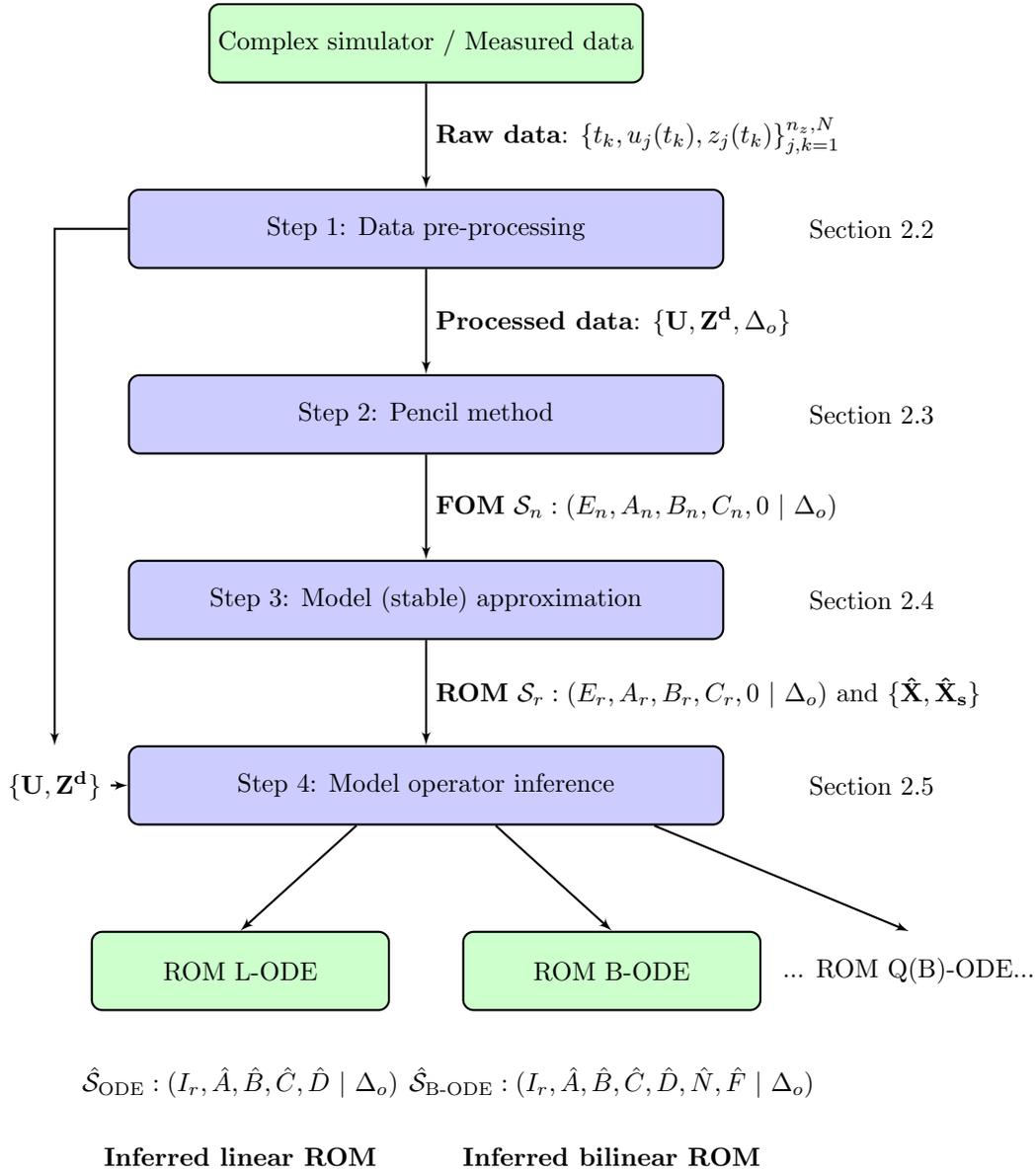
\begin{figure}
    \centering
    \tikzstyle{block} = [draw, thick, fill=blue!20, rectangle,minimum height=3em, minimum width=8cm, rounded corners]
\tikzstyle{inout} = [draw, thick, fill=green!20, rectangle, minimum height=3em, minimum width=3cm,rounded corners]
\tikzstyle{outout} = [draw, thick, fill=green!20, rectangle, minimum height=3em, minimum width=4cm,rounded corners]
\tikzstyle{block2} = [fill=white!20, rectangle,minimum height=3em, minimum width=1cm, rounded corners]
\tikzstyle{donnees} = [coordinate]
\tikzstyle{pinstyle} = [pin edge={to-,thick,black}]
\tikzstyle{connector} = [->,thick]

\begin{tikzpicture}[auto, node distance=2.5cm,>=latex']
    \node [inout] (input) {Complex simulator / Measured data};
    \node [block, below of=input] (step0) {Step 1: Data pre-processing};
    \node [block, below of=step0] (step1) {Step 2: Pencil method};
    \node [block, below of=step1] (step2) {Step 3: Model (stable) approximation};
    \node [block, below of=step2] (step3) {Step 4: Model operator inference};
    \node [block2, left of=step3, node distance=5cm] (outLeftStep3) {$\{\mathbf U,\mathbf{Z^d}\}$};
    \node [outout, below of=step3, right of=step3,node distance=2.5cm] (outRight) {ROM B-ODE};
    \node [outout, below of=step3, left of=step3,node distance=2.5cm] (outLeft) {ROM L-ODE};
    \node [block2, right of=outRight,node distance=4cm] (outLeft2) {... ROM Q(B)-ODE...};
    \node [block2, below of=outRight, node distance=1.5cm] (outRightText) {$\Hrealr_{\text{B-ODE}}:(I_r,\Ar,\Br,\Cr,\Dr,\Nr,\Fr ~\vert~ \Delta_o)$};
    \node [block2, below of=outRightText, node distance=1cm] (outRightText) {\textbf{Inferred bilinear ROM}};
    \node [block2, below of=outLeft, node distance=1.5cm] (outLeftText) {$\Hrealr_{\text{ODE}}:(I_r,\Ar,\Br,\Cr,\Dr ~\vert~ \Delta_o)$};
    \node [block2, below of=outLeftText, node distance=1cm] (outLeftText2) {\textbf{Inferred linear ROM}};
    
    \node [block2, right of=step0, node distance=6cm] (sec1) {Section \ref{ssec:data}};
    \node [block2, right of=step1, node distance=6cm] (sec2) {Section \ref{ssec:pencil}};
    \node [block2, right of=step2, node distance=6cm] (sec3) {Section \ref{ssec:rom}};
    \node [block2, right of=step3, node distance=6cm] (sec4) {Section \ref{ssec:infer}};

    \draw [connector] (input.south) -- node {\textbf{Raw data}: $\{t_k,u_{j}(t_k),z_{j}(t_k)\}_{j,k=1}^{n_z,N}$} (step0.north);
    \draw [connector] (step0.south) -- node {\textbf{Processed data}: $\{\mathbf U,\mathbf{Z^d},\Delta_o\}$} (step1.north);
    \draw [connector] (step1.south) -- node {\textbf{FOM} $\Hreal_n:(\E_n,\A_n,\B_n,\C_n,0 ~\vert~ \Delta_o)$} (step2.north);
    \draw [connector] (step2.south) -- node {\textbf{ROM} $\Hreal_r:(\E_r,\A_r,\B_r,\C_r,0 ~\vert~ \Delta_o)$ and $\{\mathbf{\hat X},\mathbf{\hat X_s}\}$} (step3.north);
    \draw [connector] (step3.-150) -- node {} (outLeft.north);
    \draw [connector] (step3.-30) -- node {} (outRight.north);
    \draw [connector] (step3.-10) -- node {} (outLeft2.north);
    \draw [connector] (outLeftStep3.east) -- node {} (step3.west);
    \draw [connector] (step0.west) -| node {} (outLeftStep3.north);
    
\end{tikzpicture}
    \caption{Proposed non-intrusive reduced order model construction process. From the time-domain input and output data collected on the complex simulator or a rich measurement set, one obtains the \textbf{ROM} linear, bilinear, quadratic.. models.}
    \label{fig:process}
\end{figure}

To construct the simplified linear or nonlinear dynamical model, a non-intrusive time-domain data-driven approach is considered here. The  main steps of the rationale are summed-up below and detailed in the rest of the section (Figure \ref{fig:process} provides an illustration of it).
\begin{enumerate}
    \item Collect the constantly sampled time-domain input and output raw data from the simulator and extract some physical information (section \ref{ssec:data}).
    \item Construct a linear dynamical \textbf{FOM} (or generating dynamical model) that reproduces as accurately as possible the time-domain raw data by catching the main linear dynamics only (section \ref{ssec:pencil}). 
    \item Compute a linear dynamical \textbf{ROM} that approximates the \textbf{FOM} one using any data or model-based interpolatory approach, leading to a simplified input to states to output mapping function (section \ref{ssec:rom}).
    \item Optimise the \textbf{ROM} by inferring an adjusted \textbf{L-ROM} or \textbf{B-ROM} (or any linear and nonlinear structure) through a mixed collection of the raw data and the reduced states obtained from the \textbf{ROM} to account for possible nonlinearities (section \ref{ssec:infer}).
\end{enumerate}
The flow described in Figure \ref{fig:process} is detailed in the rest of the section. As each step would require a dedicated attention, specific details are skipped and let to the reader curiosity through appropriate references. The main contribution of the proposed scheme stands in proposing a complete approach to deal with reduced order model inferring from limited number of time-domain input-output data. Moreover, as shown along the rest of the section, each step is scalable and leads to a process adapted to simulators embedding very complex dynamics, which also stands as a contribution for practitioners and end users.

\subsection{Collecting data and estimate delays}
\label{ssec:data}

From the time-domain simulation, let us collect at each time index $k$ the constantly sampled signals with period  $h$ ($h\in\Real_+$) input $\u_k$ and output $\z_k$  defined as
\begin{equation}
    \begin{array}{rcl}
    \z_k &:=& \left[z_1(t_k), z_2(t_k),\dots , z_j(t_k), \dots , z_{n_z}(t_k) \right]^T\in\Real^{n_z\times 1}\\
    \u_k &:=& \left[u_1(t_k), u_2(t_k),\dots , u_j(t_k), \dots , u_{n_u}(t_k) \right]^T\in\Real^{n_u\times 1},
    \end{array}
\end{equation}
where $u_j(t_k)$ and $z_j(t_k)$ denote the the $j$-th coordinate input and output at time $t_k$ ($k=1,\dots,N$, $N\in\mathbb Z$). 
The considered input and output raw data matrices then read
\begin{equation}
    \mathbf U :=
    \left[
    \begin{array}{cccc}
        \lvert & \lvert & & \lvert\\
        \u_1 & \u_2& \dots & \u_N\\
        \lvert & \lvert & & \lvert\\
    \end{array}
    \right]
    \in \Real^{n_u\times N}
    \text{ and }
    \mathbf Z := 
    \left[
    \begin{array}{cccc}
        \lvert & \lvert & & \lvert\\
        \z_{1} & \z_{2} & \dots & \z_{N}\\
        \lvert & \lvert & & \lvert\\
    \end{array}
    \right]
    \in \Real^{n_z\times N}.
\label{eq:UZ}
\end{equation}

\begin{remark}[About the delay estimation]
In multiple physical dynamical systems such as transport equations, input-output delays naturally appear. This phenomena is really specific and would be important to evaluate and to embed in the inferred model one is looking for (rather than approximating it)\footnote{Note that the same kind of comment holds true for polynomial terms considered in \cite{BennerInfer:2020}.}. In this case we suggest to compute the number of sampling delays $\tau_j\in\mathbb Z$ as ($\varepsilon>0$)
\begin{equation}
    \begin{array}{rcl}
    \tau_{j} &:=& \displaystyle\arg \max_{k} \norm{z_{j}(t_k)}\\
    && \text{s.t. } \norm{z_{j}(t_k)} < \varepsilon
    \end{array} .
\end{equation}
In other words, the transport delay for each signal sequence $z_{j}(t_k)$ is estimated. Consequently, we can define the shifted "raw" time-vector as, 
\begin{equation}
\begin{array}{rcl}
    z^d_{j} &:=& \left[z_{j}(t_{k+\tau_{j}}),z_{j}(t_{k+1+\tau_{j}}),\dots,z_{j}(t_{N-\tau_{j}}) \mathbf 1_{1\times \tau_{j}}\right]^T \in \Real^{N \times 1}\\
    \z^d_k &:=& \left[z^d_{1}(t_k), z^d_{2}(t_k),\dots , z^d_{j}(t_k), \dots , z^d_{n}(t_k) \right]^T\in\Real^{n\times 1},
\end{array}
\end{equation}
where $\mathbf 1_{1\times \tau_{j}}$ is a 1-column vector of length $\tau_{j}$ used to fill the dimension of the vector $z^d_{j}$ so that it has the same dimension $n_z$ as $z_j$. Each $\z^d_k$ vector is then the same signal as $\z_k$ but shifted to remove the delay part. The shifted "raw" output measurement matrix is now defined as
\begin{equation}
    \mathbf{Z^d} := 
    \left[
    \begin{array}{cccc}
        \lvert & \lvert & & \lvert\\
        \z^d_{1} & \z^d_{2} & \dots & \z^d_{N}\\
        \lvert & \lvert & & \lvert\\
    \end{array}
    \right]
    \in \Real^{n_z\times N},
    \label{eq:Zd}
\end{equation}
and we can define the output delay operator $\Delta_o(\cdot)$ as
\begin{equation}
\begin{array}{rccl}
    \Delta_o : & \Real^{n_z} &\rightarrow&\Real^{n_z} \\
               & \z(t_k) &\rightarrow& \z(t_k-\tau) = \z^d(t_k),
\end{array}
\label{eq:Delta}
\end{equation}
where $\tau=[\tau_1,\tau_2,\dots,\tau_{j},\dots, \tau_{n_z}]^T\in\mathbb Z^{n_z\times 1}$. This step is optional and in the general case one can consider $\mathbf{Z^d}:=\mathbf Z$ and $\tau:=\mathbf{0}_{n_z \times 1}$. 
\end{remark}

\subsection{Full order model construction}
\label{ssec:pencil}


Following the second step of the process illustrated in Figure \ref{fig:process}, on the basis of the collected input data $\mathbf{U} \in \Real^{n_u\times N}$, shifted output data $\mathbf{Z^d} \in \Real^{n_z\times N}$ and output delay operator $\Delta_o$, we now want to construct a $n$-th order \mimo \textbf{FOM} linear dynamical function $\Htran\in\Cplx^{n_z\times n_u}$, defined as 
\begin{equation}
    \tilde \z = \Htran\u,
\end{equation}
which maps inputs $\u$ to the outputs $\tilde \z$. One objective is that feeding $\Htran$ with the time-domain data collected in $\mathbf U$ results in a time-series $\tilde \z(t_k)$ as close as possible to the raw data gathered in $\mathbf{Z}$. Here we seek for a transfer function $\Htran$ described by its realisation $\Hreal_n:(\E_n,\A_n,\B_n,\C_n,0 ~\vert~ \Delta_o)$ explained as
\begin{equation}
    \Hreal_{n} : \left \lbrace
    \begin{array}{rcl}
    \E_{n} \x_{n}(t_k+h) &=& \A_{n} \x_{n}(t_k) + \B_{n} \u(t_k)\\
    \tilde \z(t_k) &=&\Delta_o(\C_{n} \x_{n}(t_k))
    \end{array}
    \right. ,
    \label{eq:real_fom}
\end{equation}
where $\x_n(t_k)\in\Real^{n}$, $\u(t_k)\in\Real^{n_u}$ and $\tilde \z(t_k)\in\Real^{n_z}$ are the (global) internal, input and output variables gathering all the measured outputs\footnote{In the rest of the paper we consider $n_u=1$ for simplicity. The rest of the results remain true.}. The $\E_n$, $\A_n$, $\B_n$ and $\C_n$ matrices are real and of appropriate dimension, resulting from the sub-generating models $\Hreal_j$ described hereafter. In the proposed scheme, the computation of these matrices may be complex when the dimensions $n_z$, $n_u$ and $N$ (and thus $n$) increase, which is the  considered case in this work. This problem is generally referred to as a \emph{model  identification} one for which a complete research field is attached (see \eg  \cite{Schoukens:2016} for a very complete survey). One greedy but yet effective way to construct the generating \textbf{FOM} given in \eqref{eq:real_fom} is by identifying the $n_z$ \miso transfers independently. Then, the \textbf{FOM} matrices have the following structure,
\begin{equation}
\begin{array}{rcl}
\E_n &=& \textbf{blkdiag}(\E_1,\E_2,\dots,\E_{n_z})\\
\A_n &=& \textbf{blkdiag}(\A_1,\A_2,\dots,\A_{n_z})\\
\B_n &=& [\B_1^T,\B_2^T,\dots,\B^T_{n_z}]^T\\
\C_n &=& \textbf{blkdiag}(\C_1,\C_2,\dots,\C_{n_z})
\end{array}.
\label{eq:structure}
\end{equation}
where each $j$-th sub matrix has a dimension induced by the Pencil methods.  One then seeks for $\Htran_j\in\Cplx^{1\times n_u}$, the sub-generating \miso model linking the input $\u$ to each output $\tilde z_j$ ($j=1,\dots,n_z$). Here we consider that each transfer $\Htran_j$  with realisation  $\Hreal_{j}:(\E_{j},\A_{j},\B_{j},\C_{j},0 ~\vert~ \tau_{j})$ is 
\begin{equation}
    \Hreal_{j} : \left \lbrace
    \begin{array}{rcl}
    \E_{j} \x_{j}(t_k+h) &=& \A_{j} \x_{j}(t_k) + \B_{j} \u(t_k)\\
    \tilde z_{j}(t_k) &=&\C_{j} \x_{j}(t_k-\tau_{j})
    \end{array}
    \right. ,
    \label{eq:real_fom_j}
\end{equation}
where $\x_{j}(t_k)\in\Real^{n_{j}}$, $\u(t_k)\in\Real^{n_u}$ and $\tilde z_{j}(t_k)\in\Real$ are (local) internal, input and output variables related to the transfer from $\u$ to $\tilde z_j$, the $j$-th measurement point. Being given the delay $\tau_j$, the construction of the rest of the realisation is performed by using the $j$-th line of the shifted data $\mathbf{Z^d}$ (the delay term being added afterward). Here, the \emph{pencil method} is employed to construct \eqref{eq:real_fom_j}. This method relies on a specific re-arrangement of the shifted raw output data $\mathbf{Z^d}$ followed by any rank revealing factorisation. Interested reader is invited to refer to the original paper \cite{HoKalman:1966} or more recent comprehensive tutorials \cite{IonitaPencil:2012,AntoulasPLOSOne:2018} for details on the constructions steps. One of the main property of the associated $n_j$-th transfer function  $\Htran_{j}(z)=\C_{n}(z\E_{n}-\A_{n})^{-1}\B_{n}z^{-\tau_j}$ is that its impulse response generates the sequence $\tilde z_{j}(t_k)$ that approximates the $j$-th line of the raw vector $\mathbf{Z}$. Moreover, the order $n_j$ is encoded in the ($\A_j,\E_j$) pencil and revealed by the \textbf{SVD} factorisation.

Then by stacking the $j=1,\dots,n_z$ sub-generating \miso models $\Hreal_{j}$ linking each input to each output, one obtains \eqref{eq:real_fom}. The resulting dimension is
\begin{equation}
    n = \displaystyle\sum_{j=1}^{n_z}n_{j}.
    \label{eq:dim}
\end{equation}

\begin{remark}[About $\Hreal_n$ complexity]
At this point, $\Hreal_n$ embeds important complexity  mostly due to the number of considered grid points increasing the matrices dimension as explained in \eqref{eq:dim}. This complexity prevents any linear algebra operation such as eigenvalue computation, impulse response computation... Still, as the $\E_n$ and $\A_n$ \eqref{eq:structure} elements are mostly sparse, one can exploit this property for dimensionality reduction. This is illustrated later on in the application section \ref{sec:appli}.
\end{remark}

\subsection{Stable reduced order model construction}
\label{ssec:rom}

Now, one has access to a dynamical model that well approximates the main linear dynamics of the simulator, but which also embeds a important complexity. Before considering the possible discarded nonlinear terms (in step 4 - next section), this \textbf{FOM} is first simplified by a \textbf{ROM} to alleviate computational limitations (step 3 - this section). 

\subsubsection{Settings and context}

Following the third step of Figure \ref{fig:process}, rooted on the \textbf{FOM} $\Hreal_n$ given in \eqref{eq:real_fom}, which dimension increases with the number of measurement grid points $n_z$, one now desires to compute $\Hreal_r:(\E_r,\A_r,\B_r,\C_r,0 ~\vert~ \Delta_o)$, the reduced order model given as
\begin{equation}
    \Hreal_r :\left \lbrace
    \begin{array}{rcl}
    \E_r  \x_r(t_k+h) &=& \A_r \x_r(t_k) + \B_r \u(t_k)\\
    \overline \z(t_k) &=&\Delta_o(\C_r \x_r(t_k))
    \end{array}
    \right. .
    \label{eq:real_all_red}
\end{equation}
where $\x_r(t_k)\in\Real^{r}$, $\u(t_k)\in\Real^{n_u}$ and $\overline \z(t_k)\in\Real^{n_z}$ are the reduced internal, input and approximated output variables. This problem is the so-called \emph{model reduction} one. The goal generally is to approximate the original model with a smaller and simpler one, having the same structure and similar response characteristics as the original. For an overview of model reduction methods, we refer the reader to the book \cite{AntoulasBook:2020}. Without entering into detailed considerations, two broad reduction families exist. The model-driven (generally projection-based) and the data-driven ones. Both usually consider minimising the mismatch error between $\Hreal_n$ and $\Hreal_r$ using some frequency-domain (or complex-domain) criteria. To avoid too much complexity and notations we invite reader to refer to \eg \cite{AntoulasBook:2020,Peherstorfer:2017,GugercinSIAM:2008,VanDooren:2008,Mayo:2007} or to references given in section \ref{sec:intro} for more details.  

In the considered use-case, under mild technical and practical considerations, we will use the data-driven Loewner interpolatory framework, originally presented in \cite{Mayo:2007} (see also a tutorial in \cite{AntoulasSurvey:2016}). Main lines are recalled hereafter for completeness.

\subsubsection{The Loewner interpolation}

The Loewner framework is a \emph{data-driven} method aimed at building a $m$-th order rational descriptor \lti dynamical model $\Htran_m$ which interpolates given complex data. In the proposed scheme, the previously generated \textbf{FOM} $\Htran_n$ without delay, \ie with realisation $(\E_n,\A_n,\B_n,\C_n,0  ~\vert~ \mathbf{0}_{n_z\times 1})$ (or shorthand $(\E_n,\A_n,\B_n,\C_n,0)$) is considered. Let the left (or row) data be given together with the right (or column) data, as below
\begin{equation}
\left.
\begin{array}{c}
(\mu_j,\lv_j^H,\vv_j^H) \\
\text{for $j=1,\dots,m$}
\end{array}
\right\}
\text{~~and~~}
\left\{
\begin{array}{c}
(\lambda_i,\rv_i,\wv_i) \\
\text{for $i=1,\dots ,m$}
\end{array}
\right. ,
\label{eq:loewnerInput}
\end{equation}
where $\vv_j^H=\lv_j^H\Htran_n(\mu_j)$ and $\wv_i=\Htran_n(\lambda_i)\rv_i$, with $\lv_j\in\Cplx^{n_y\times 1}$, $\rv_i\in\Cplx^{n_u\times 1}$, $\vv_j\in\Cplx^{n_u\times 1}$ and $\wv_i\in\Cplx^{n_y\times 1}$. In addition, the set of distinct interpolation points $\{ z_k\}_{k=1}^{2m} \subset \Cplx$ is split up into two equal subsets ($\lambda_i, \ \mu_j\in\Cplx$), \ie
\begin{equation}\label{eq:shift}
\{z_k\}_{k=1}^{2m} =\{\mu_j\}_{j=1}^{m} \cup \{\lambda_i\}_{i=1}^{m}.
\end{equation}
The method then consists in building the \emph{Loewner} matrix $\LL \in \Cplx^{m\times m}$ and \emph{shifted Loewner} matrix $\sLL \in \Cplx^{m\times m}$ defined as follows, for $i=1,\dots,m$ and $j=1,\dots,m$:
\begin{align}\label{eq:loewnerMatrices}
\begin{split}
[\LL]_{j,i} &= \dfrac{\vv_j^H\rv_i - \lv_j^H\wv_i}{\mu_j - \lambda_i} 
= \dfrac{\lv_j^H\big( \Htran_n(\mu_j) - \Htran_n(\lambda_i) \big) \rv_i}{\mu_j - \lambda_i}, \\
\,[\sLL]_{j,i} &= \dfrac{\mu_j\vv_j^H\rv_i - \lambda_i\lv_j^H\wv_i}{\mu_j - \lambda_i}
= \dfrac{ \lv_j^H\big( \mu_j\Htran_n(\mu_j) - \lambda_i\Htran_n(\lambda_i) \big) \rv_i}{\mu_j - \lambda_i}.
\end{split}
\end{align}
Then, similarly to the pencil method, the model $\Htran_m$ given by the descriptor realisation,
\begin{equation}
\Hreal_m:\left \lbrace
\begin{array}{rcl}
E_m  \x(t_k+1) &=& A_m \x(t_k) + B_m \u(t_k)\\
\y(t_k) &=&C_m \x(t_k)
\end{array}
\right. ,
\label{eq:loewnerDescrR}
\end{equation}
where $E_m = -\LL$, $A_m = -\sLL$, $[B_m]_k = \vv_k^H$ and $[C_m]_k = \wv_k$ (for $k=1,\ldots,m)$, with the related transfer function $\Htran_m(z) = C_m(z E_m-A_m)^{-1}B_m$, interpolates $\Htran_n$ at the given driving frequencies and directions defined in \eqref{eq:loewnerInput}, \ie satisfies the conditions
\begin{equation}
\begin{array}{rcl}
\lv_j^H\Htran_m(\mu_j) &=& \lv_j^H\Htran_n(\mu_j) \\
\Htran_m(\lambda_i)\rv_i &= &\Htran_n(\lambda_i) \rv_i
\end{array}.
\label{eq:loewnerIntep}
\end{equation}
Assuming that the number $2m$ of available data is large enough and that (for $k=1,\ldots,2m$)
\begin{equation}
    \rank (z_k \LL - \sLL) = 
    \rank ([\LL,\sLL]) = 
    \rank ([\LL^H,\sLL^H]^H) = n_-,
    \label{eq:rankCond}
\end{equation}
where $z_k$ are as in \eqref{eq:shift}, then it is shown in \cite{Mayo:2007} that a minimal order model $\Htran_{n_-}$ of dimension $n_- \leq m$ that still satisfies the interpolatory conditions \eqref{eq:loewnerIntep} can be computed by projecting \eqref{eq:loewnerDescrR}. Let $Y \in \Cplx^{m \times n_-}$ be the matrix containing the first $n_-$ left singular vectors of $[\LL,\sLL]$ and $X \in \Cplx^{m \times n_-}$  the matrix containing the first $n_-$ right singular vectors of $[\LL^H,\sLL^H]^H$, then,
\begin{equation}
\label{eq:proj}
    \E_{n_-} = Y^H \E_m X,\,
    \A_{n_-} = Y^H \A_m X,\, 
    \B_{n_-} = Y^H \B_m \text{ and }
    \C_{n_-} = \C_m X,
\end{equation}
is a realisation $\Hreal_{n_-}$ of the $\Htran_{n_-}(z) = \C_{n_-}( z \E_{n_-}-\A_{n_-})^{-1}\B_{n_-}$, encoding a \emph{minimal McMillan degree} equal to $\rank(\LL)$. The quadruple $\Hreal_{n_-}:(E_{n_-},A_{n_-},B_{n_-},C_{n_-},0)$ is a descriptor realization of $\Htran_{n_-}$. Note that if ${n_-}$ in \eqref{eq:rankCond} is greater than $\rank(\LL)$, then $\Htran_{n_-}$ can either have a direct-feedthrough term or a polynomial part. Finally, the number ${n_-}$ of singular vectors composing $Y$ and $X$ used to project the system $\Htran_{n_-}$ in \eqref{eq:proj} may be decreased to $r<{n_-}$ at the cost of deteriorating the interpolatory conditions, leading to the reduced model $r$-th order rational model denoted $\Htran_r(z)= \C_r(z\hat \E_r-\A_r)^{-1}\B_r$ with realisation $\Hreal_r:(\E_r,\A_r,\B_r,\C_r,0)$. This allows a trade-off between complexity of the resulting model and accuracy of the interpolation.

\begin{remark}[About the output delay operator $\Delta_o(\cdot)$]
Even if the above model can be constructed including the output delay, the \textbf{ROM} realisation $\Hreal_r:(\E_r,\A_r,\B_r,\C_r,0)$ is computed on the basis of the delay-free \textbf{FOM} data. Therefore, to recover the model \eqref{eq:real_all_red}, the output delay operator $\Delta(\cdot)$ should be added afterward. 
\end{remark}

\subsubsection{Comments model input-output stability}

While the Loewner framework ensures some interpolatory properties with a minimal realisation (observable and controllable), it does not provide any guarantee concerning the stability of the resulting matrix pencil $(\E_r,\A_r)$. In the perspective of time-domain (predictive) simulation, this is a major issue. Therefore, a post treatment should be applied to approximate the Loewner-based interpolant by a stable one. To this aim, the recent work of \cite{Kohler:2014} is used to project a rational unstable model onto its best stable approximant. It consists in projecting the rational model $\Htran_r$ onto its closest stable subset, here using the $\Hinf$-norm, leading to a sub-optimal stable model of the same dimension. Mathematically, given a realisation $\Hreal_r$ associated to $\Htran_r\in \Linf$, one aims at finding $P_\infty(\Htran_r) \in \Hinf$ such that,
\begin{equation}
P_\infty(\Htran_r) = \arg \inf_{\mathbf G \in \Hinf}\norm{\Htran_r-\mathbf G}_{\mathcal H_\infty}.
\label{eq:stableApprox}
\end{equation}
Technical details and assumption can be found in \cite{Kohler:2014}. The new projected model $\Htran_r\leftarrow P_\infty(\Htran_r)$ is now stable and close to the unstable original one. Its impulse response should be close to the original data, but accuracy losses may be expected for some outputs.

\subsection{Reduced order model inference}
\label{ssec:infer}

Now one has access to a \textbf{ROM}  $\Hreal_r$ being  representative of the \textbf{FOM}  $\Hreal_n$. However, as already pointed in the construction of  $\Hreal_n$, the \textbf{ROM}  may capture the linear dynamics only and poorly the nonlinear ones. Moreover, the reduction and post-stablisation steps may have been done at the cost of a deteriorated accuracy. This is why now one aims  at inferring a (linear or nonlinear) model by taking advantage of this  \textbf{ROM} model, as a complement of the raw data \eqref{eq:UZ}. This last step will first adjust the model, and if necessary, suggest the addition of a nonlinear term.

\subsubsection{Settings and context}

At this point, by simply simulating the impulse response of the stable $\Hreal_r$ model, one has access to a reduced state-space data set in addition to the original raw data. This last point is the key enabler to apply operator inference as in \cite{Peherstorfer:2015}. Following the fourth step of Figure \ref{fig:process} and by assuming $\E_r$ is invertible (which was already the case in the previous stable projection step  \cite{Kohler:2014}), the reduced model $\Hreal_r$ can be re-written as
\begin{equation}
    \Hreal_r :\left \lbrace
    \begin{array}{rcl}
    \xr(t_k+h) &=& \A_r \xr(t_k) + \B_r \u(t_k)\\
    \overline \z(t_k) &=&\Delta_o(\C_r \xr(t_k)+ \D_r \u(t_k))
    \end{array}
    \right. .
\end{equation}
where $\A_r\leftarrow \E_r^{-1}\A_r$ and $\B_r\leftarrow \E_r^{-1}\B_r$ and potentially non-null $\D_r\in\Real^{n_z\times n_u}$ term resulting from the $\E_r$ matrix rank (details are skipped here but reader may refer to \cite{AntoulasSurvey:2016}). The internal state is now denoted as $\xr(t_k)\in\Real^r$. Now the $r$-th order \textbf{ROM} is available, we are able to easily simulate its outputs in response to the raw input data $\mathbf U$ given in \eqref{eq:UZ}. We can now obtain the approximated outputs
\begin{equation}
\mathbf{\overline Z} := 
    \left[
    \begin{array}{cccc}
        \lvert & \lvert & & \lvert\\
        \overline\z_{1} & \overline\z_{2} & \dots & \overline\z_{N}\\
        \lvert & \lvert & & \lvert\\
    \end{array}
    \right]
    \in \Real^{n_z\times N}.
\end{equation}
and, more interestingly, the reduced state and shifted state trajectories as
\begin{equation}
\begin{array}{rcl}
    \mathbf{\hat{X}} &:= &
    \left[
    \begin{array}{cccc}
        \lvert & \lvert & & \lvert\\
        \xr_{1} & \xr_{2} & \dots & \xr_{N-1}\\
        \lvert & \lvert & & \lvert\\
    \end{array}
    \right]
    \in \Real^{r\times (N-1)} \text{ and } \\
    \mathbf{\hat{X}_s} &:= &
    \left[
    \begin{array}{cccc}
        \lvert & \lvert & & \lvert\\
        \xr_{2} & \xr_{3} & \dots & \xr_{N}\\
        \lvert & \lvert & & \lvert\\
    \end{array}
    \right]
    \in \Real^{r\times (N-1)}.
\end{array}
\label{eq:stateTraj}
\end{equation}

These reduced state data \eqref{eq:stateTraj} are of low dimension and will be the starting point for inferring a model using the input data $\mathbf U$, shifted raw data $\mathbf{Z^d}$ and reduced state ones $\mathbf{\hat{X}}$ and $\mathbf{\hat{X}_s}$. 

\begin{remark}[About classical operator inference]
In the seminal paper \cite{Peherstorfer:2015}, authors present an approach for operator inference. In  \cite{Peherstorfer:2015}, the authors assume that one has access to the full state snapshots matrix from the original simulator. Then, authors apply a \textbf{POD} to project the original states onto a reduced subspace. They also assume to be able evaluating / post-processing these state snapshot thought a known function to recover the non-linearity. Still, as an important result,  authors also demonstrate the convergence between the  intrusive (projection-based) and non-intrusive resulting operators, by only assuming stabilised snapshots and convergent numerical scheme. Here instead, the proposed approach assumes that the input $\mathbf U$ and measurements $\mathbf Z$ snapshots are the only known variables. The states considered here result from the above constructed \textbf{ROM}. This difference makes the problem somehow in between operator inference and identification results. 
\end{remark}

In the following, the operator inference step is synthesized following the same assumptions as in  \cite{GoseaDMD:2020} and \cite{Peherstorfer:2015}. The convergence and proofs are also let to reader curiosity.

\subsubsection{Structured linear model inference}

Now we have access to the reduced state $\mathbf{\hat{X}}$ and shifted states $\mathbf{\hat{X}_s}$ trajectories \eqref{eq:stateTraj}, together with the delay shifted raw data $\mathbf{Z^d}$, one can infer a model by solving a least-square problem. Following \cite{GoseaDMD:2020} one can solve the following problem
\begin{equation}
    \min_{\Ar,\Br,\Cr,\Dr} \norml{
    \vectortwo{\mathbf{{\hat{X}_s}}}{\mathbf{Z^d}} - \matrixtwo{\Ar}{\Br}{\Cr}{\Dr} \vectortwoT{\mathbf{\hat{X}}}{\mathbf U}
    }_F,
    \label{eq:ls_lODE}
\end{equation}
to recover the $r$-th order reduced matrices $(\Ar,\Br,\Cr,\Dr)$. Note that solving the above least square problem may probably modify the eigenvalues of $\Ar$ and may lead (once again) to a loss of stability. This is why the following structured problem  \eqref{eq:ls_lODEstruct} may be considered instead (quite similar with the \textbf{DMD} with control of \cite{Proctor:2016})
\begin{equation}
    \min_{\Br,\Cr,\Dr} \norml{
    \vectortwo{\mathbf{{\hat{X}_s}}-\Ar\mathbf{\hat{X}}}{\mathbf{Z^d}} - \matrixtwo{0_{r\times r}}{\Br}{\Cr}{\Dr} \vectortwoT{\mathbf{\hat{X}}}{\mathbf U}
    }_F, 
    \label{eq:ls_lODEstruct}
\end{equation}
where $\Ar$ is fixed and typically chosen as $\Ar=\A_r$. In this second case, the dynamical matrix $\Ar$ is considered as known and input/output residual adjustment is applied only. Obviously, it would result in a lower accuracy but will preserve the stability. Solving \eqref{eq:ls_lODEstruct} may be viewed as a residual adjustment step since the dynamical matrix is left unchanged. This final step, after adding the output delay, then leads to the inferred reduced linear ordinary difference equations \textbf{L-ODE} given by the following realisation $\Hrealr_{\textbf{L-ODE}}:(I_r,\Ar,\Br,\Cr,\Dr  ~\vert~ \Delta_o)$
\begin{equation}
    \Hrealr_{\textbf{L-ODE}} :\left \lbrace
    \begin{array}{rcl}
    \xr(t_k+h) &=& \Ar \xr(t_k) + \Br \u(t_k)\\
    \zr(t_k) &=&\Delta_o(\Cr \xr(t_k) + \Dr\u(t_k))
    \end{array}
    \right. ,
    \label{eq:lODE}
\end{equation}
which is both stable and of reduced dimension. This latter may then be used in place of the original simulator for fast evaluations. 

\subsubsection{Structured nonlinear model inference}

So far, all the considered models were linear. It is also possible at this stage to enrich the model structure to account for the variability and the specificity of the phenomena and to deal with potential nonlinear contributions that may not be caught in the previous steps. This can be done \eg through the addition of a bilinear term. Then one solves \eqref{eq:ls_bODE},
\begin{equation}
    \min_{\Ar,\Br,\Cr,\Dr,\Nr,\Fr} \norml{
    \vectortwo{\mathbf{{\hat{X}_s}}}{\mathbf{Z^d}} - 
    \left[\begin{array}{ccc}
    \Ar & \Br & \Nr \\
    \Cr & \Dr & \Fr 
    \end{array}\right]
    \left[\begin{array}{ccc} 
    \mathbf{\hat{X}} & \mathbf{U} & \mathbf{\hat{X}U}
    \end{array}\right]
    }_F, 
    \label{eq:ls_bODE}
\end{equation}
where $\Nr\in\Real^{r\times r}$ and $\Fr\in\Real^{n_z\times r}$. Interestingly, enriching the model structure with the bilinear term allows enhancing the model restitution. Note that in the literature, discrete-time model with this form are also called affine systems (see very complete work from \cite{BoydPhD:1985}). The interest of this affine structure is that, provided the nonlinear operator generating the data has a fading memory, then it can be approximated to arbitrary accuracy with a model of this form (see chapter 4 of \cite{BoydPhD:1985}). In that case, its impulse response can be approximated by a Volterra series which is precisely the response of affine models. Here again for stability issues, instead of solving \eqref{eq:ls_bODE}, one may similarly consider solving the problem \eqref{eq:ls_bODEstruct}, putting the nonlinearity on the output equation only, \ie setting $\Nr=0_{r\times r}$ and $\Ar=\A_r$,
\begin{equation}
    \min_{\Br,\Cr,\Dr,\Fr} \norml{
    \vectortwo{\mathbf{{\hat{X}_s}}-\Ar\mathbf{\hat{X}}-\Nr\mathbf{\hat{X}}\mathbf U}{\mathbf{Z^d}} - 
    \left[\begin{array}{ccc}
    0_{r\times r} & \Br & 0_{r\times r} \\
    \Cr & \Dr & \Fr 
    \end{array}\right]
    \left[\begin{array}{ccc} 
    \mathbf{\hat{X}} & \mathbf{U} & \mathbf{\hat{X}U}
    \end{array}\right]
    }_F.
    \label{eq:ls_bODEstruct}
\end{equation}

By adding the output delay, this leads to the so-called affine or bilinear ordinary difference equations \textbf{B-ODE} given by the following realisation $\Hrealr_{\textbf{B-ODE}}:(I_r,\Ar,\Br,\Cr,\Dr,0_{r\times r},\Fr ~\vert~ \Delta_o)$
\begin{equation}
    \Hrealr_{\textbf{B-ODE}} :\left \lbrace
    \begin{array}{rcl}
    \xr(t_k+h) &=& \Ar \xr(t_k) + \Br \u(t_k)\\
    \zr(t_k) &=&\Delta_o(\Cr \xr(t_k) + \Dr\u(t_k)+ \Fr\xr(t_k)\u(t_k))
    \end{array}
    \right. .
    \label{eq:bODEstruct}
\end{equation}

At this point, the additional bilinear term in the output equation of \eqref{eq:bODEstruct} is used to enhance the output signal restitution. Similarly to the remark related to problems \eqref{eq:ls_lODE} and  \eqref{eq:ls_lODEstruct}, one can also seek for a bilinear term in the state equation but with the risk of obtaining an unstable model. This constraints depends on the available experimental data. If the raw data converge to an equilibrium, this constrain may be unnecessary. However, if the simulation has not necessarily converged, this latter should be considered. Of course, extensions may consider specific forms for $\Ar$ and $\Nr$ (\eg Schur) to address this stability issue in a constrained optimisation step. Still, this problem remains out of the scope of the proposed contribution and would be considered in future works. To close this model inferring section, one should notice that the problems \eqref{eq:ls_lODE}, \eqref{eq:ls_lODEstruct}, \eqref{eq:ls_bODE} and \eqref{eq:ls_bODEstruct} can all be solved using the \textbf{DMD} approach developed in  \cite{GoseaDMD:2020}. And, as a straightforward extension one may also consider quadratic and/or bilinear models (see \cite{GoseaDMD:2020}).

Now the complete process introduced in Figure \ref{fig:process} has been detailed, the next section \ref{sec:appli} is attached to apply it on a very complex pollutants simulator use-case.

\section{Application: pollutant dispersion analysis and prediction}
\label{sec:appli}
Air pollution is a growing concern in many countries around the world. According to the World Health Organization (WHO), 90~\% of the population was exposed to pollutants concentrations transgressing the WHO guidelines in 2016 \cite{WHO2016}. Efficient modeling tools are therefore essential to forecast air pollution episodes and apply the appropriate emissions control strategies for limiting their intensity. 
The main challenge consists in having reliable models with reasonable computational costs to be used in  forecast mode. In this context, simplified dynamical models as the one presented in section \ref{sec:inference} could represent very suitable tools. 
A brief introduction about atmospheric pollutants dispersion is presented in section \ref{sec:intro:context}. 
Then, the application of the proposed \textbf{ROM} modeling presented in section \ref{sec:inference} is applied on an atmospheric pollution use-case. The data are collected from a \textbf{LES} presented in section \ref{sec:mnh}. To illustrate the process, we first consider a rough grid in section \ref{ssec:process}, focusing on the four steps of the scheme given in  Figure \ref{fig:process}. Then, in section \ref{ssec:simulation} we present simulation results over a thin grid, illustrating the accuracy of the inferred nonlinear \textbf{ROM} on this very complex problem.

\subsection{Context of pollutants dispersion}
\label{sec:intro:context}

The pollutants dispersion is driven by multiple processes ranging over several spatio-temporal scales \cite{mayer1999air}. These processes can be classified into three categories :
\begin{itemize}
    \item the spatio-temporal distribution of the sources of contaminants which can be released by human activities routinely  (industry, road traffic, residential) or accidentally (explosion, fire);
    \item the meteorological situation which includes turbulent dispersion, convective mixing, transport by the wind; 
    \item the chemical reactions of the atmospheric compounds and the photo-chemistry. 
\end{itemize}

A detailed representation of all the processes is needed to ensure an accurate forecast of the plume behavior at fine scale. 
Several kinds of tools can be used to simulate pollutants dispersion depending on the the spatio-temporal scales involved. 
On the one hand, Gaussian or Lagrangian models allow to represent long term and annual trends at a low computational cost through their simple parametrisations of the atmospheric boundary layer \cite{unal2005airport,peace2006identifying}. 
On the other hand, meso-scale atmospheric models include a more elaborated representation of the physical processes thereby allowing a more accurate modeling of the plume dispersion at fine scale \cite{pison2004quantification,woody2015estimates}.  Besides, these models can be run under light wind conditions -- which represent the most favorable situations for pollutant accumulation -- while Gaussian or Lagrangian approximations do not apply in these cases.
Although atmospheric models appear more appropriate for an accurate representation of the plume dynamics, their implementation is very onerous as the computational cost increases as the resolution becomes finer. 
In this context, the development of new methods for reducing these models may represent a major improvement as it would provide fast and accurate responses to extreme cases of pollution exposition. 
To our knowledge, this is the first time that a dynamical model is constructed in order to reproduce the plume dispersion obtained from a \textbf{LES}. The development of such a model is challenging as the plume behavior is driven by several nonlinear processes. A description of the \textbf{LES} configuration is proposed in the next section.

\subsection{Dedicated simulator description}
\label{sec:mnh}

The \textbf{LES} is run with the non-hydrostatic atmospheric research model \textbf{Meso-NH} \cite{lac2018overview} which has already been successfully used to model air quality and dispersion  \cite{Sarrat:2017,sarrat2006impact,lac2013co, auguste2020large}. \textbf{Meso-NH} resolves the unsteady 3D Euler equations under the anelastic hypothesis. A wide range of parametrisation is available to represent the processes of microphysics, chemistry, aerosols, radiation, convection, turbulence and surface-atmosphere interactions. Curious readers may refer to the scientific documentation for a complete overview of the model abilities (\texttt{http://mesonh.aero.obs-mip.fr/mesonh54/BooksAndGuides}).  \\

\begin{table}
\centering
\begin{tabular}{ccp{0.1mm}}
\toprule
Domain size & N$_X$ = 900  \hspace{4mm}  N$_Y$ = 900  \hspace{4mm}  N$_Z$  = 79 & \\     
Horizontal resolution & $\delta$x = $\delta$y = 10~m &\\
Vertical resolution  & $\delta$z = 2~m near the ground &\\
\midrule
Lateral boundary conditions & \multicolumn{1}{c}{\begin{tabular}[c]{@{}c@{}} Cyclic for meteorological variables \\ Open for passive tracers \end{tabular}} & \\\cmidrule{2-3} 
Upper boundary condition & Rayleigh absorbing layer \cite{lafore1997meso}  & \\
\midrule
Atmosphere initialization & \multicolumn{1}{c}{\begin{tabular}[c]{@{}c@{}} Potential temperature and humidity profiles \\ Westerly logarithmic wind profile  \end{tabular}} & \\
\midrule
Surface scheme & SURFEX \cite{masson2013surfexv7} & \\\cmidrule{2-3} 
Transport scheme & \multicolumn{1}{c}{\begin{tabular}[c]{@{}c@{}} CEN4TH/RKC4 for momentum \\ PPM for meteorological and scalar variables \cite{colella1984piecewise} \end{tabular}} & \\\cmidrule{2-3} 
Turbulence scheme & \multicolumn{1}{c}{\begin{tabular}[c]{@{}c@{}} 3D scheme of order 1.5 \cite{cuxart2000turbulence}  \\  Deardorff mixing length \cite{deardorff1974three}   \end{tabular}} & \\\cmidrule{2-3} 
Radiation scheme & \multicolumn{1}{c}{\begin{tabular}[c]{@{}c@{}} Fouquart and Bonnel eq. for short-waves \cite{fouquart1980computations}  \\  RRTM for long-waves \cite{mlawer1997radiative}  \end{tabular}} & \\
\bottomrule
\end{tabular}
~\\
\caption{Summary table of the numerical configuration used to perform the \textbf{LES}.}
\label{table:numconfig}
\end{table}

 \begin{figure}
    \centering
    \includegraphics[width=.45\columnwidth]{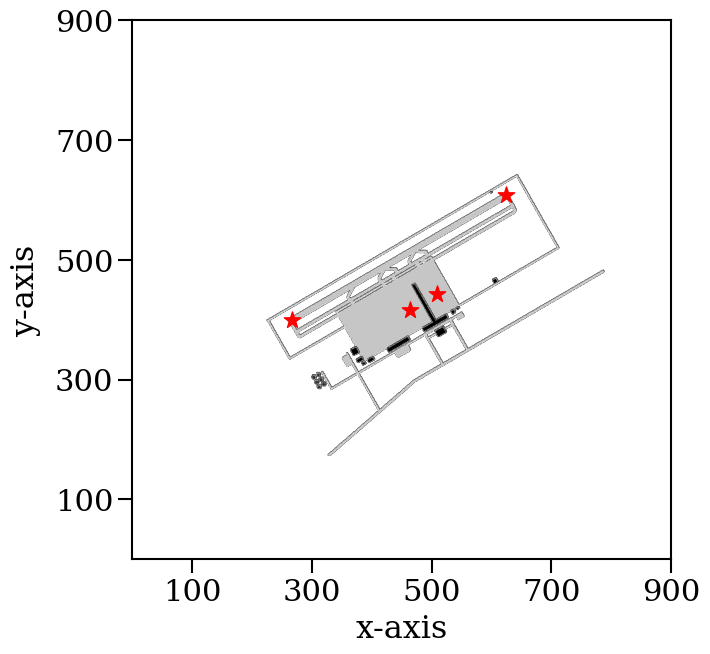}
    \caption{Simulation domain used to perform the \textbf{LES}. The red stars indicate the passive tracers sources.}
    \label{fig:domain}
\end{figure}
 
The numerical configuration used to perform the \textbf{LES} is based on \cite{sabatier2020}  who simulated the behavior of an aircraft-related plume over a 14-hour cycle. The numerical set-up is summarized in Table \ref{table:numconfig} while the simulation domain is displayed in Figure \ref{fig:domain}. This domain extends over 900 x 900 grid points with a 10 meters horizontal resolution and is centered over an airport. The vertical grid is composed of 79 levels with a 2~m resolution within the first 30~m above ground level (agl) and a gradual stretching aloft. The simulation is integrated over 3 hours with a 1-second time step. 
The atmosphere is initialized with a westerly logarithmic wind profile characterized by an average intensity of 0.5~$m.s^{-1}$ within the first 200~m agl. 
Passive tracers are released at 2~m height agl throughout the simulation  from four fixed locations materialized by red stars in Figure \ref{fig:domain}. The tracers are emitted at an arbitrary constant flux and lead to the formation of a local pollution plume. 
The plume behavior is driven by the dynamics of the atmospheric boundary layer (turbulence, flow, convective mixing) which is represented at high resolution as the simulation is performed in \textbf{LES} mode (i.e. 70-80~\% of turbulence is dynamically solved while the remaining is computed by the turbulence schema. cf. Table \ref{table:numconfig}). The surface-atmosphere interactions are also represented at fine scale through the coupling of \textbf{Meso-NH} with the surface model \textbf{SURFEX}. Momentum, heat and water fluxes are computed at each time step by considering the characteristics of each tile (texture, roughness length...). 
The outputs of the simulation are the atmospheric tracers concentrations $z_{i,j}(t_k)$ are collected at 2~m agl over each point of the domain $\{x_i\}_{i=1}^{n_x}$ and $\{y_j\}_{j=1}^{n_y}$ with a 1-minute resolution ($t_k = 0,\dots, 180$). \\ 

\begin{figure}
    \centering
    \includegraphics[width=.45\columnwidth]{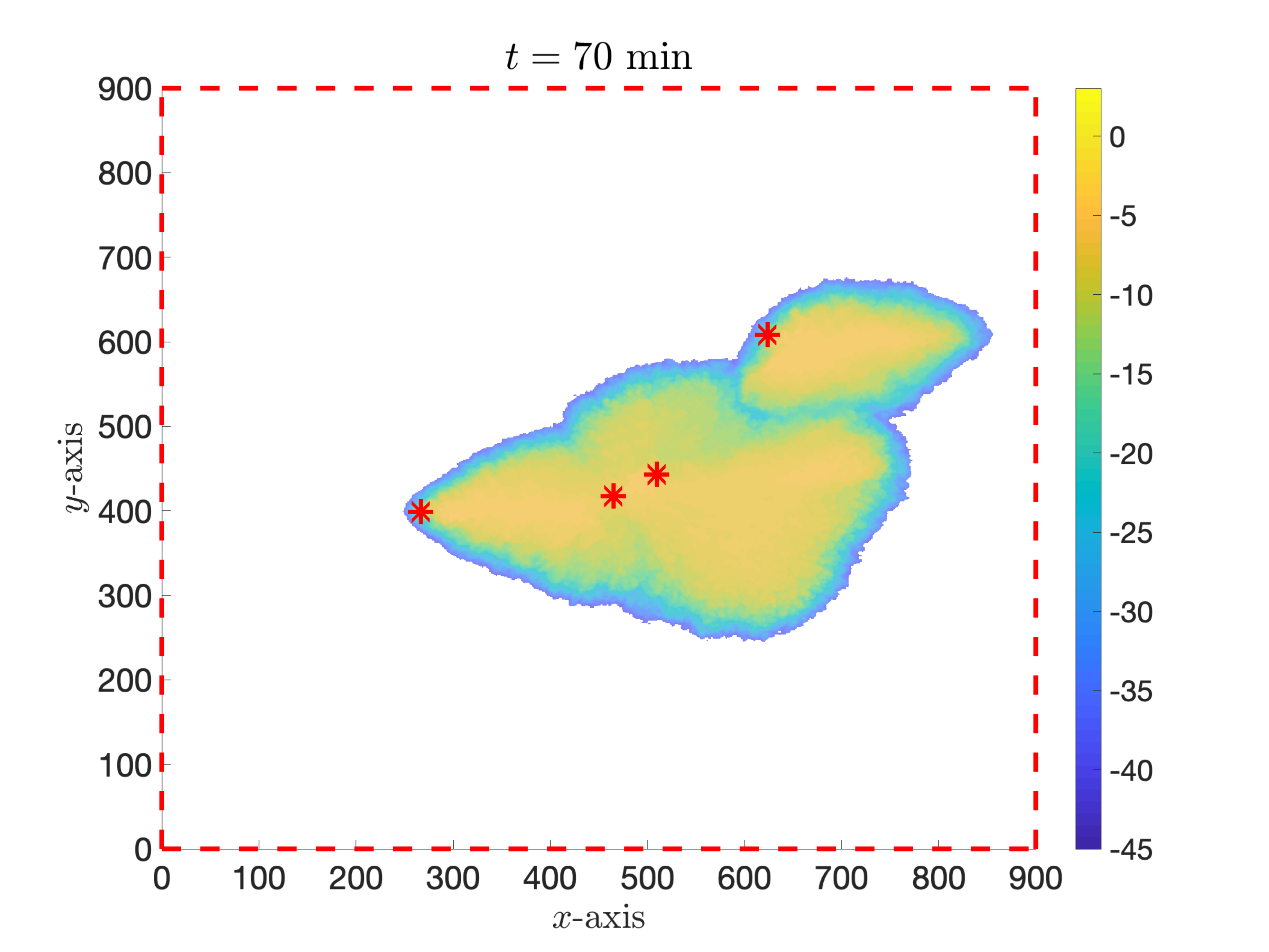}
    \includegraphics[width=.45\columnwidth]{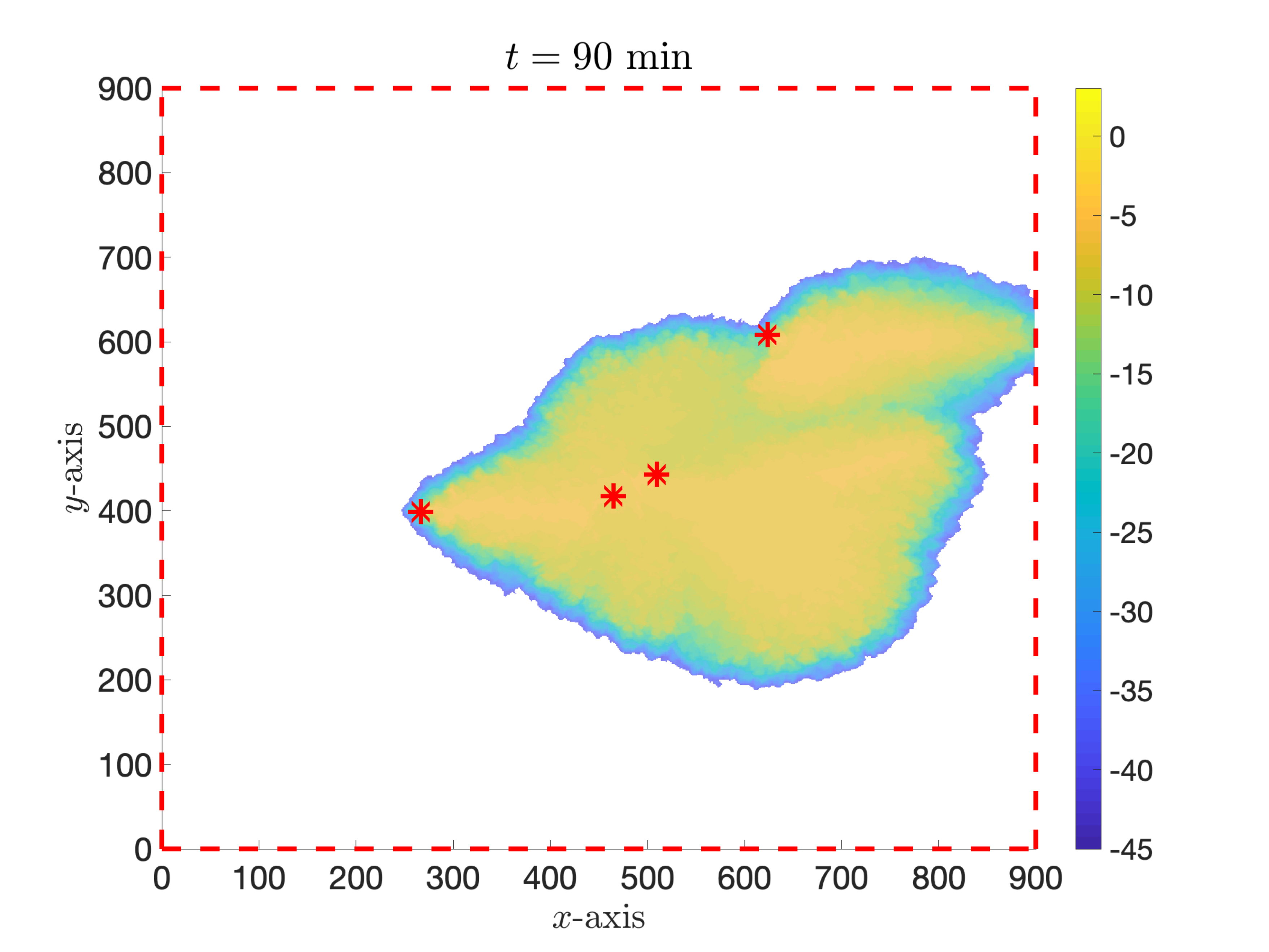}
    \includegraphics[width=.45\columnwidth]{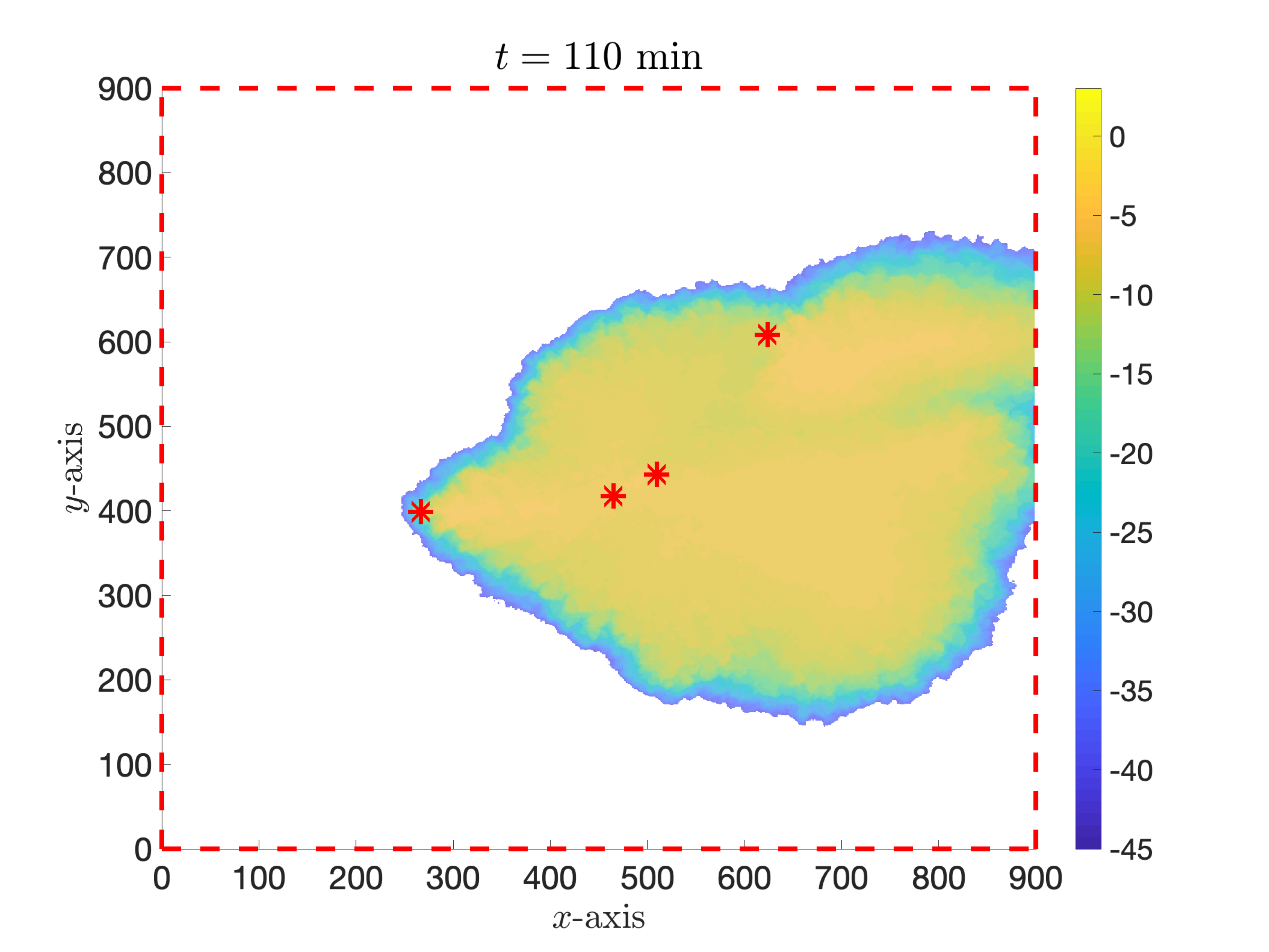}
    \includegraphics[width=.45\columnwidth]{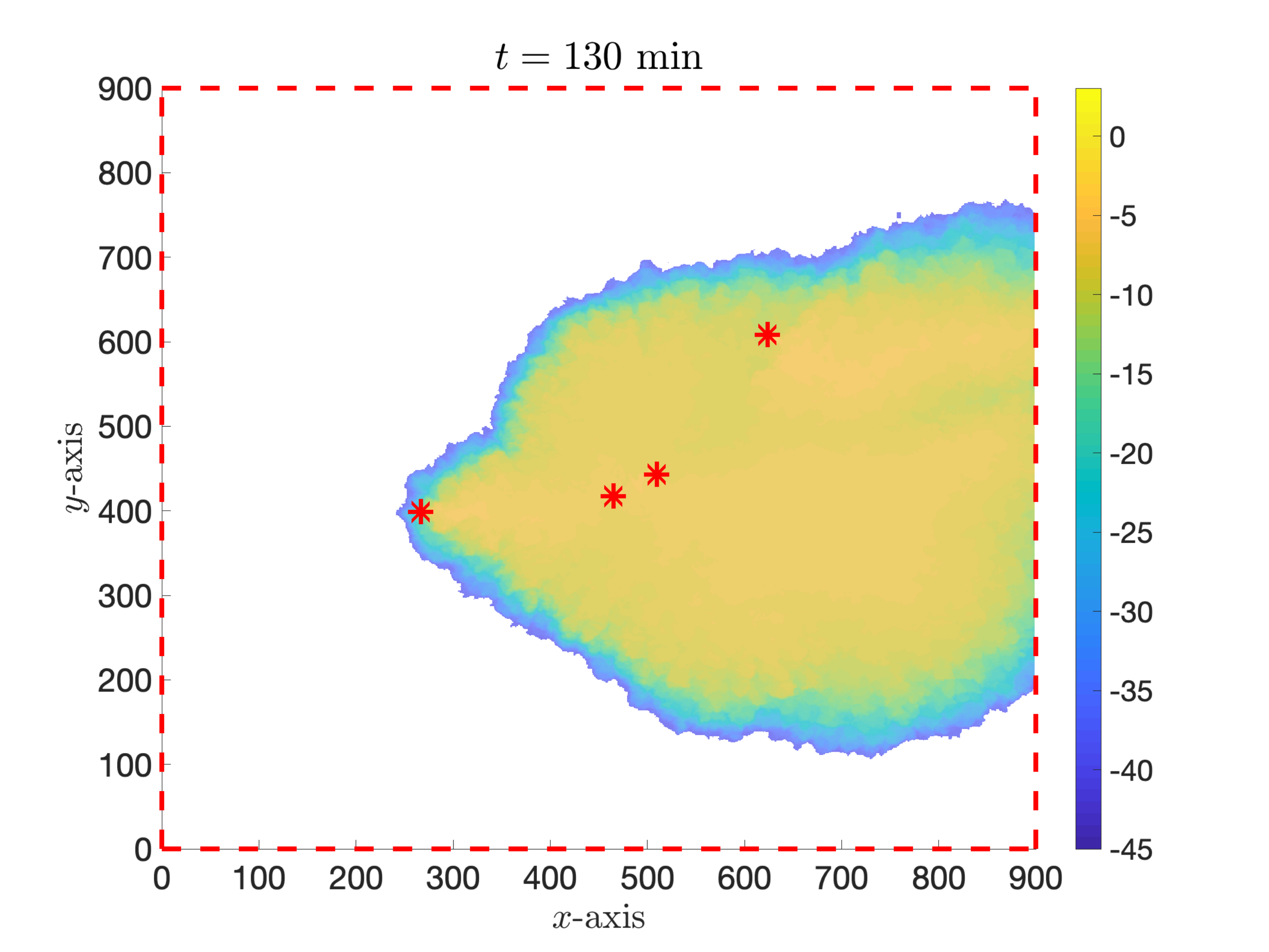}
    \includegraphics[width=.45\columnwidth]{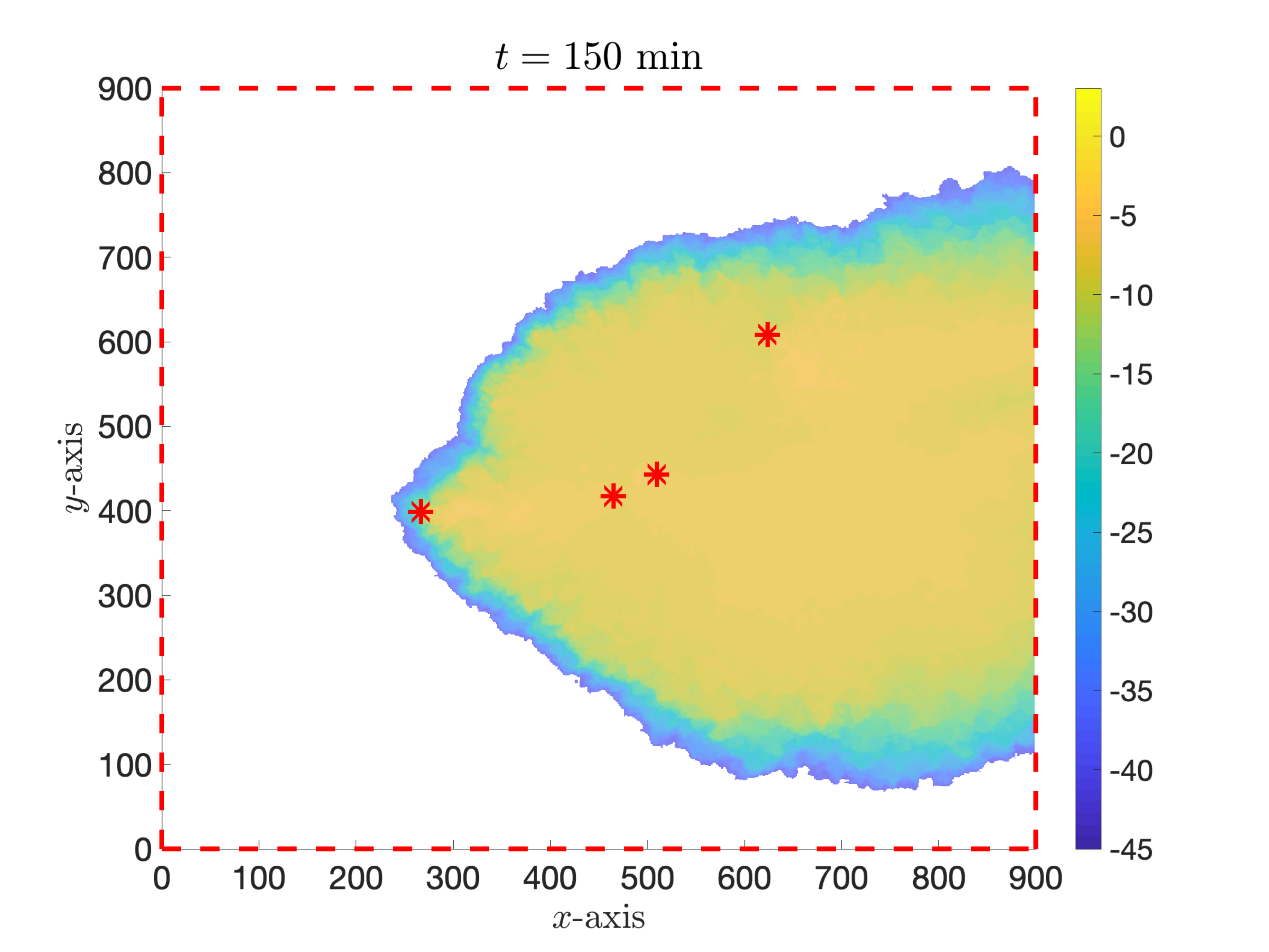}
    \includegraphics[width=.45\columnwidth]{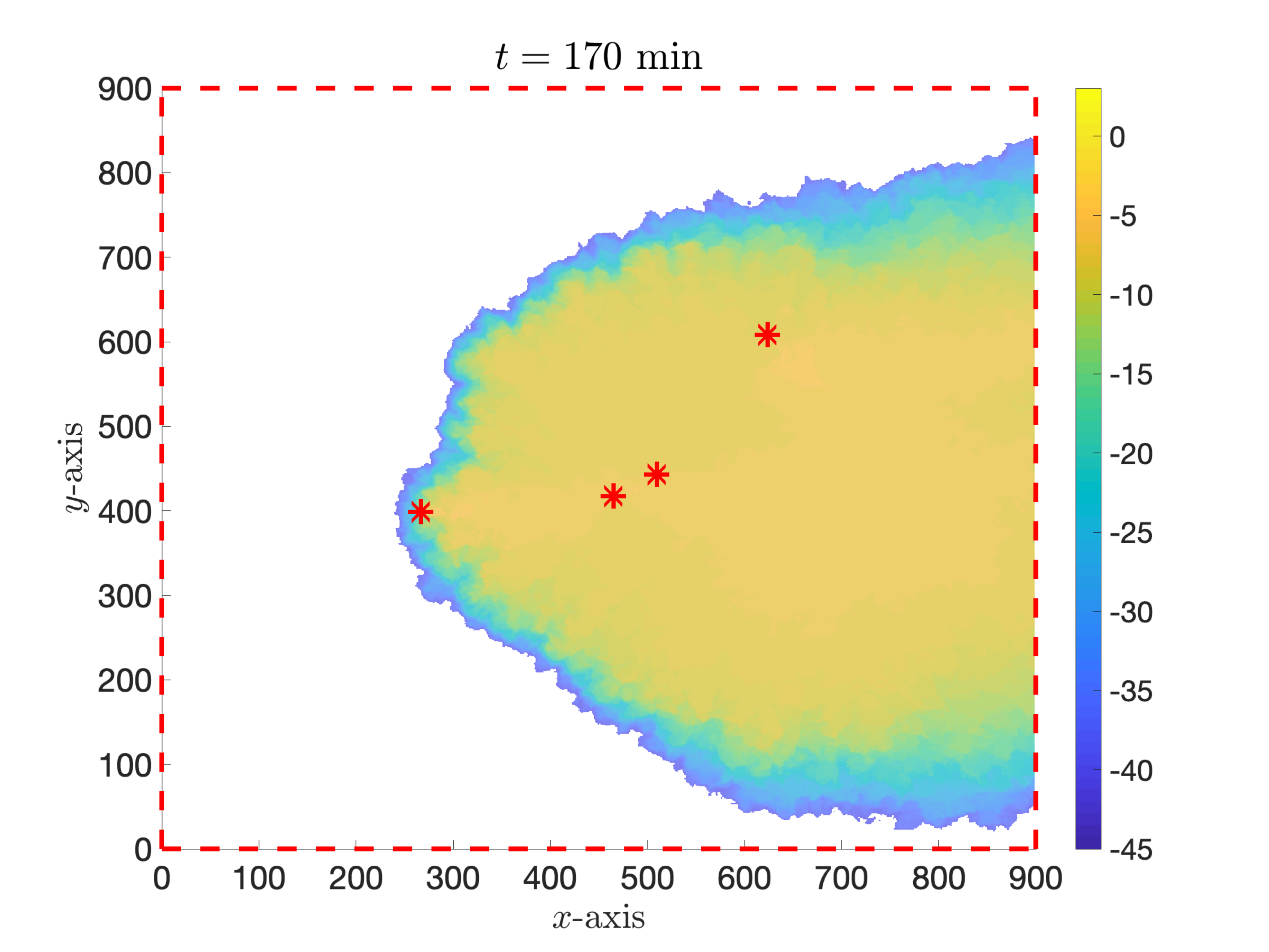}
    \caption{Snapshots of the horizontal tracers plume structure at 2~m agl from the \textbf{LES}. Pollutant values and color-bar are given with a logarithmic scale to emphasize the plume dispersion.}
    \label{fig:xyFOM}
\end{figure}

An illustration of the plume evolution is proposed in Figure \ref{fig:xyFOM} through horizontal cross-sections of tracer concentrations (expressed in a log scale) at different time instants. The plume is progressively transported eastward by the background wind while an horizontal tracers spread occurs due to the local dynamics developed at the platform scale. Regarding numerical resources, the \textbf{LES} is run over 360 processors and the 3 hours of simulation represent a total CPU cost of 5800 hours.  

\subsection{Rough grid: process illustration}
\label{ssec:process}

After running the \textbf{LES}, one collects the pollutants concentration $z_{i,j}(t_k)$ over each grid  point $\{x_i\}_{i=1}^{n_x}$ and $\{y_j\}_{j=1}^{n_y}$.  In this first part, we consider a rough $(x,y)$ spacial grid: $x_i=[600,700,800,900]$, $y_j=[300,400,500,600]$ ($n_x=n_y=4$) and a sampling-time $h=1$min with  $t_k=0,\dots,180$ ($N=180$). 

\begin{figure}
    \centering
    \includegraphics[width=.8\columnwidth]{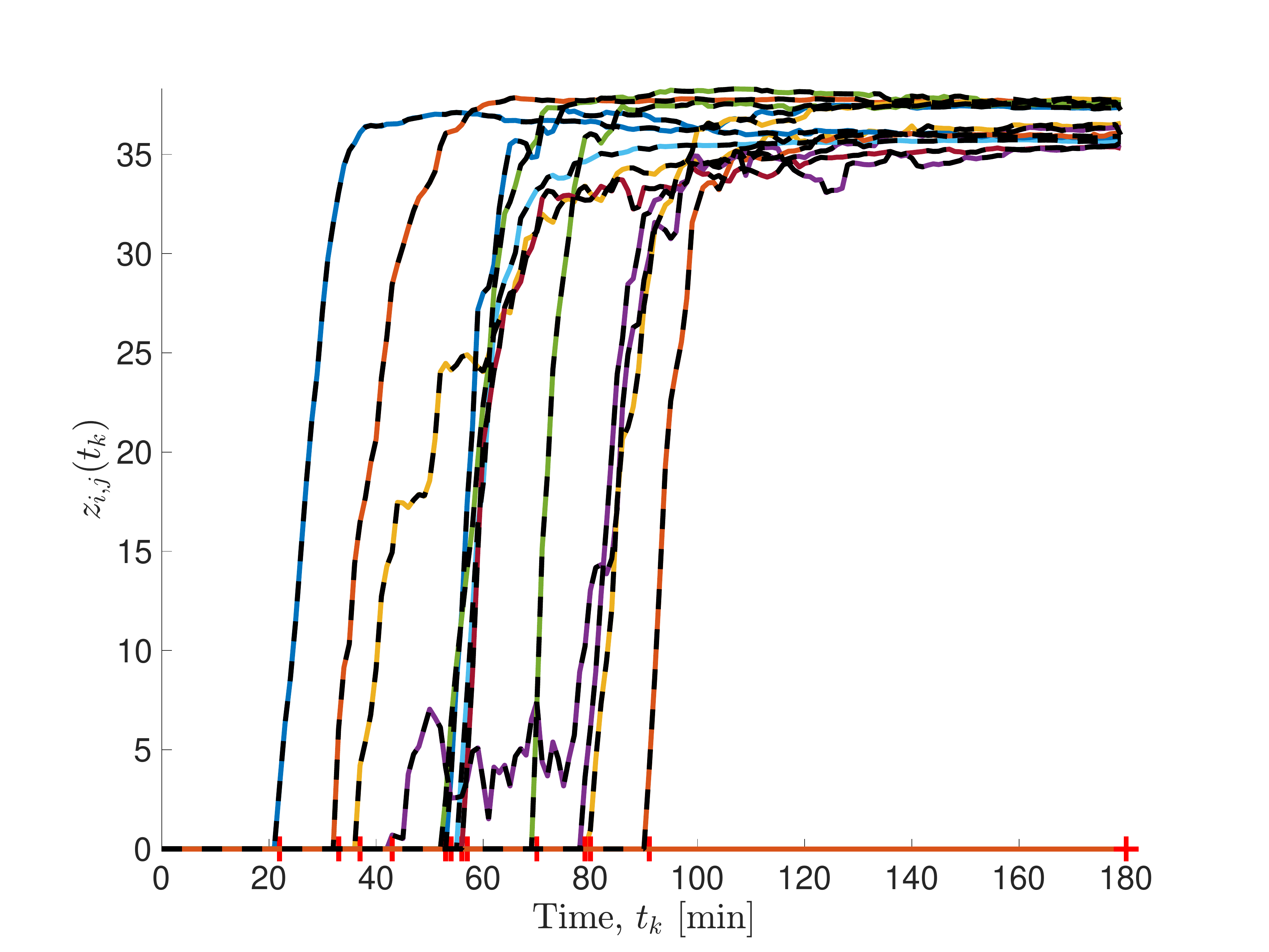}
    \caption{(Rough mesh) Original \textbf{LES} raw data $\mathbf Z$ (coloured solid lines) vs. impulse response of the \textbf{FOM} $\Hreal_n$ (black dashed). Red crosses materialise the estimated output delay values $\tau_{i,j}$. Pollutant values are given with an absolute scale starting from zero, to simplify the readability.}
    \label{fig:step1_fom}
\end{figure}

By running the first step of the process, the resulting raw input and output vectors \eqref{eq:UZ}  read $\mathbf U\in\Real^{1\times 180}$ (where $\mathbf U	\equiv 1$) and $\mathbf Z\in\Real^{16\times 180}$. Then, by estimating the sampling delay, one ends up with $\{\tau_{j}\}_{j=1}^{16}=[53,32,36,42,69,55,56,21,90,79,78,52,179,179,179,179]$min. The shifted raw data $\mathbf{Z^d}$ as in \eqref{eq:Zd} and output delay operator $\Delta_o$ as in \eqref{eq:Delta} are used then, in a second step, to construct the \textbf{FOM} $\Hreal_n$ resulting in a $n=655$-th order model as in \eqref{eq:real_fom}. Original raw data (coloured lines), estimated delay (red crosses) are reported in Figure \ref{fig:step1_fom}, and compared to the impulse response of the \textbf{FOM} (black dashed lines).

In this simplified case, the impulse response of the \textbf{FOM} alsmot perfectly recovers the original raw data. Note that the perfect match is not always observed when more output are collected exhibiting nonlinear phenomena (this will be the case in the next section). Responses may also result unstable (\ie eigenvalues of the $(\E_n,\A_n)$ pencil outside the unit circle). This point is treated in the next step. As mentioned in the above section \ref{sec:inference}, the \textbf{FOM} $\Hreal_n$ is constructed by taking all 16 sub-models $\Hreal_{j}$, leading to a sparse dynamical matrix $\A_n$, with about $5\%$ of non-zero elements (note also that the sparsity rate increases with the number $n_z$ of considered measurements).

\begin{figure}
    \centering
    \includegraphics[width=.8\columnwidth]{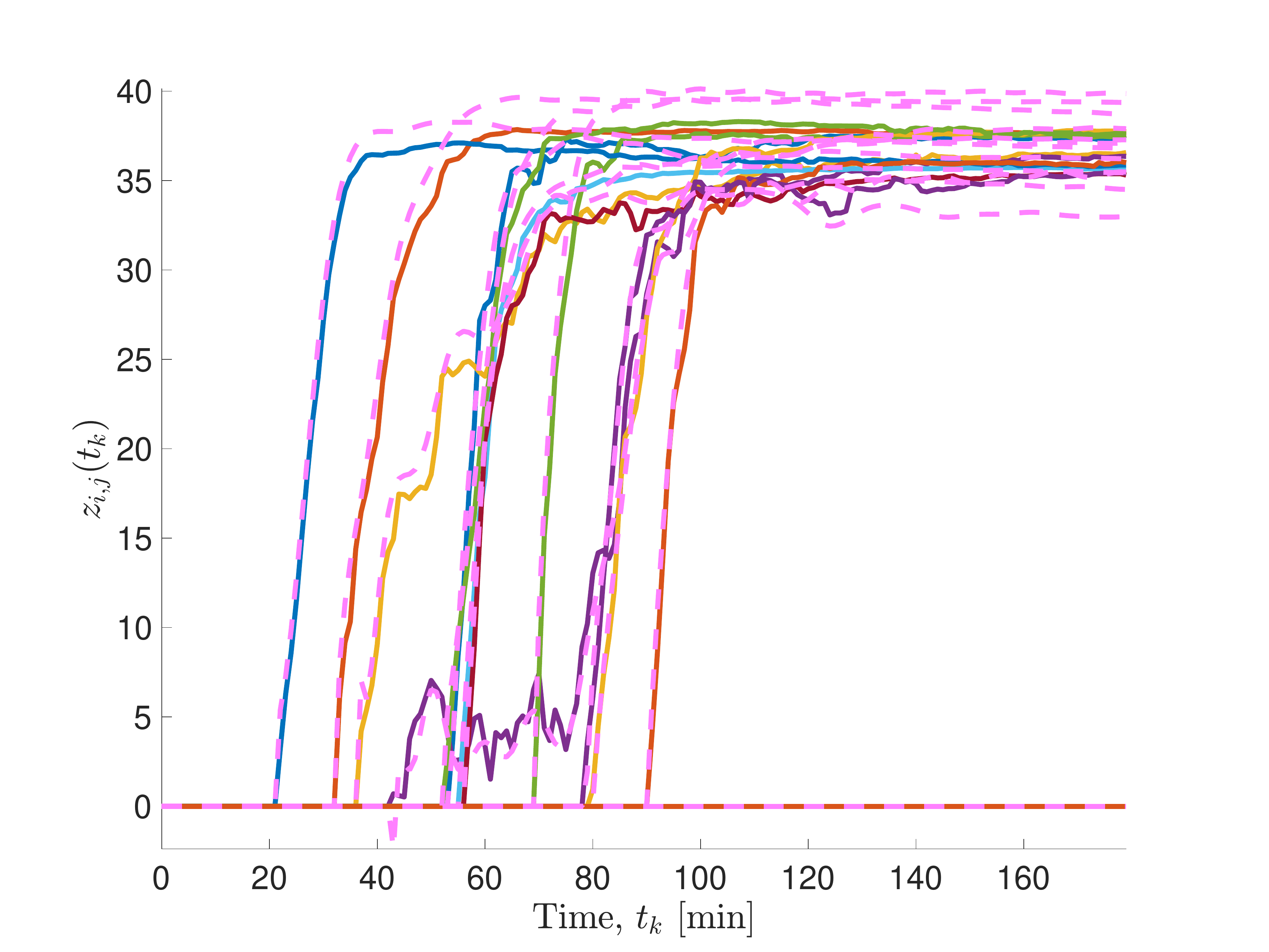}
    \caption{(Rough mesh) Original \textbf{LES} raw data $\mathbf Z$ (coloured solid lines) vs. impulse response of the stable \textbf{ROM} $\Hreal_r$ (magenta dashed). Values and color-bar are given in logarithmic scale to emphasize the plume dispersion. Pollutant values are given with an absolute scale starting from zero, to simplify the readability.}
    \label{fig:step2_rom}
\end{figure}

With only 16 measurements points, the resulting  \textbf{FOM} $\Hreal_n$ is already equipped with a realisation of dimension $n=655$. When the grid point will increase, $n$ will rapidly grow and the simulation will not be numerically feasible anymore. Thus, following the third step, a model approximation is then performed\footnote{Without entering into the details, in the considered application, the Loewner framework considers 300 sampling points (with logarithmic space) along the unit circle, from $10^{-5}$ up to $10^{-2}$.}. This step constructs a stable \textbf{ROM} $\Hreal_r$ of dimension $r=30$ (the choice of the order is subject to many considerations left apart of this paper). The resulting impulse response (dashed magenta) is compared with the raw data  in Figure \ref{fig:step2_rom}.

\begin{figure}
    \centering
    \includegraphics[width=.8\columnwidth]{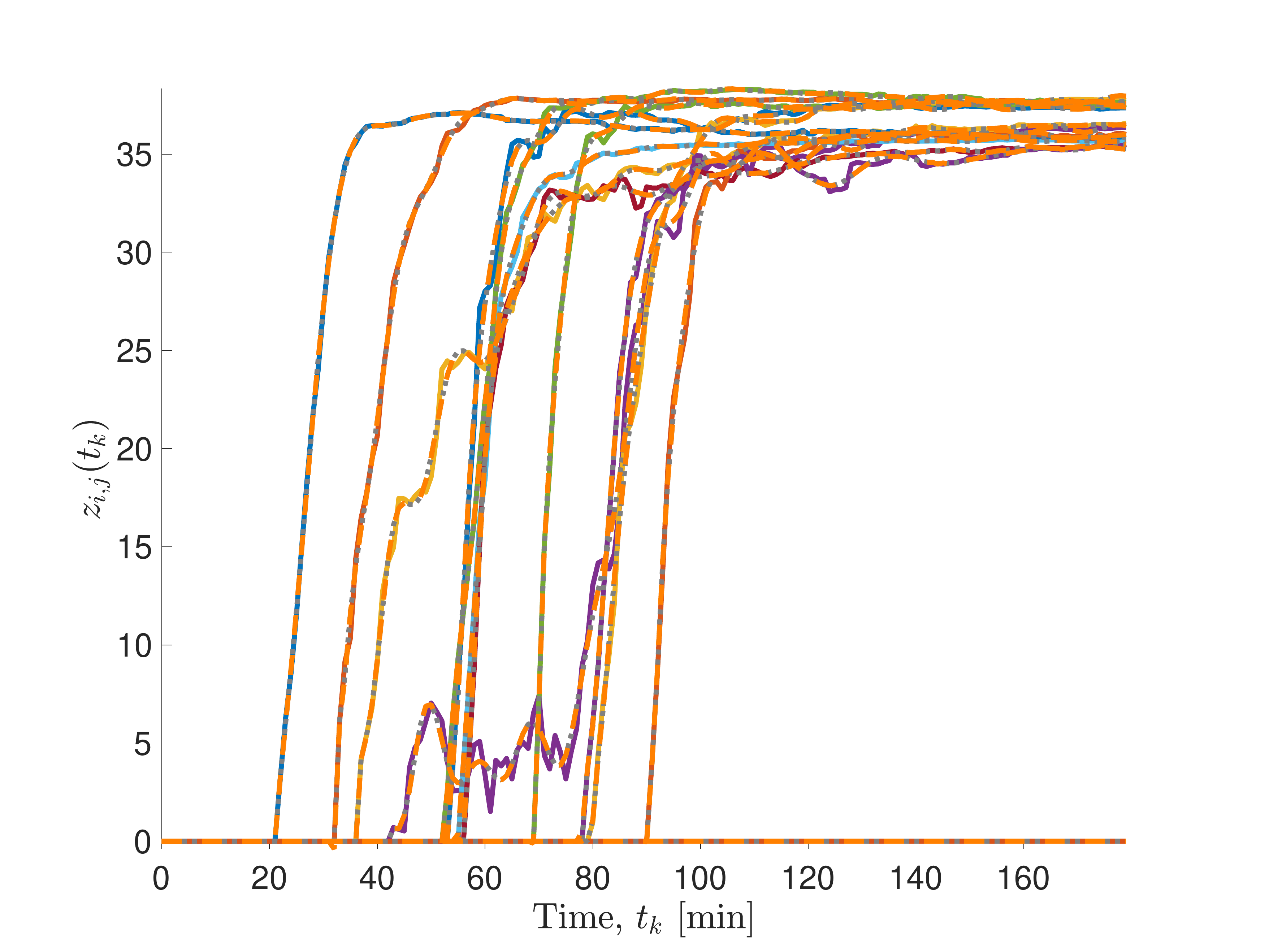}
    \caption{(Rough mesh) Original \textbf{LES} raw data $\mathbf Z$ (coloured solid lines) vs. impulse response of the inferred stable \textbf{ROM} $\Hrealr_\text{L-ODE}$ (grey dotted) and $\Hrealr_\text{B-ODE}$ (orange dashed). Pollutant values are given with an absolute scale starting from zero, to simplify the readability.}
    \label{fig:step3_inferred}
\end{figure}

\begin{figure}
    \centering
    \includegraphics[width=\columnwidth]{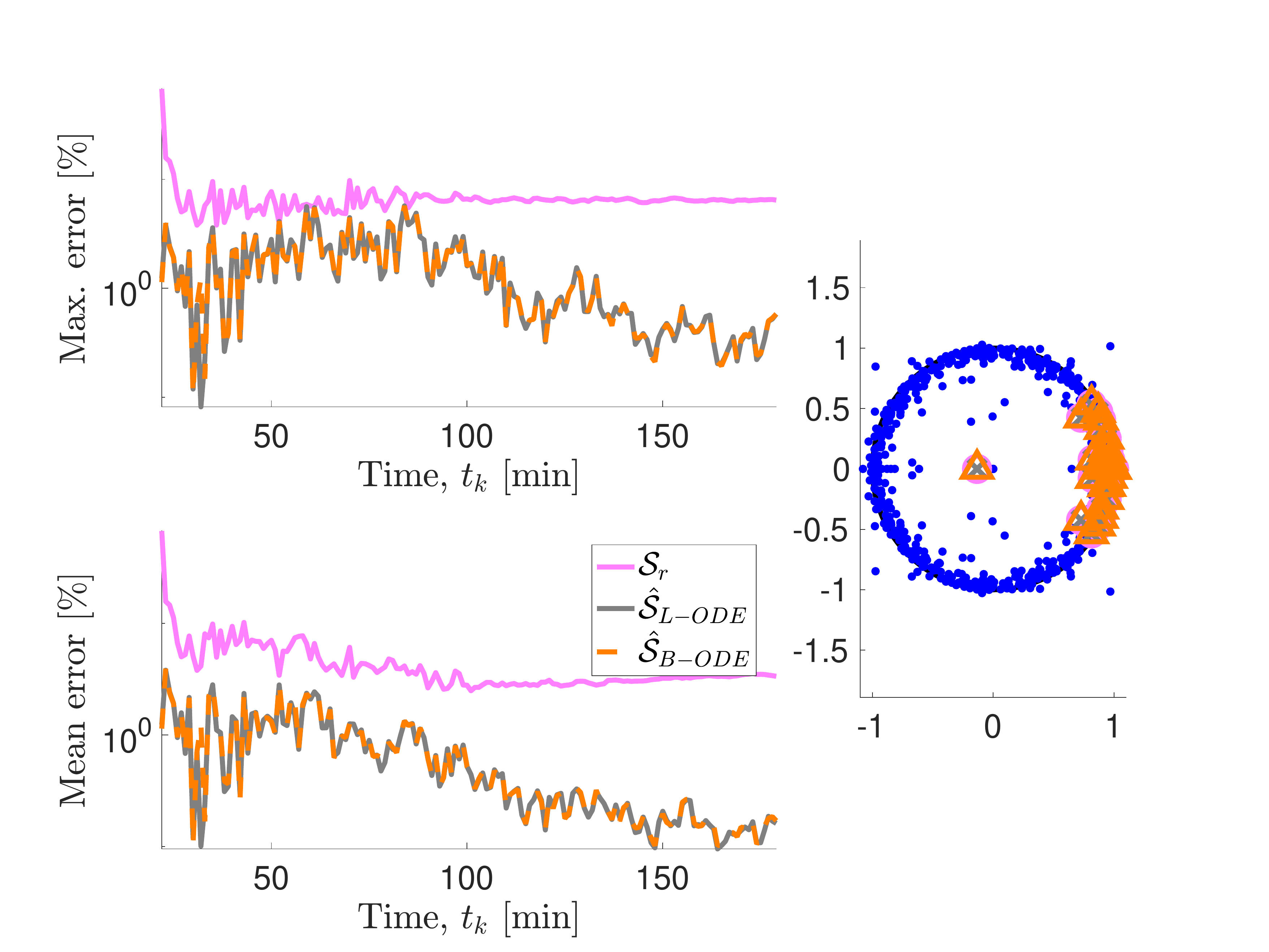}
    \caption{(Rough mesh) Top left (max) / bottom left (mean) mismatch error between the the original \textbf{LES} raw data $\mathbf Z$ and the \textbf{ROM} $\Hreal_n$ (pink), $\Hrealr_\text{L-ODE}$ (grey dotted) and $\Hrealr_\text{B-ODE}$ (orange dashed). Right: eigenvalues of $\A_n$ and $\Ar$.}
    \label{fig:step3_inferredError}
\end{figure}

Interestingly, at this point one observes that the \textbf{ROM} $\Hreal_r$ well reproduces the dynamics of the original data, but, losses in accuracy. This observation may be justified by multiple reasons: the order reduction, the post-stabilisation enforcement, the number of considered outputs, the numerical accuracy of the methods... but also to the considered criteria. Indeed, a complex-domain interpolatory approach has been used here. However, in this case, this may not be imputed to the linear structure as the \textbf{FOM} was almost exact. This latter objective may not result in a perfect restitution of all single signals but rather in a global one. Still, the Loewner framework is well adapted to very large-scale problems and is thus suited to our complex setup (see next section). At this point, the dynamics are well reproduced but not all single transfer. This is why in the fourth and final step, a least square method inspired from \cite{GoseaDMD:2020}, involving the delayed raw data $\mathbf{Z^d}$ and the reduced state-space matrices $\mathbf{\hat X}$ and $\mathbf{\hat X_s}$ generated by $\Hreal_r$, is applied. This allows  inferring a \textbf{ROM} with either a linear ($\Hrealr_{\text{L-ODE}}$) or bilinear / state-affine ($\Hrealr_{\text{B-ODE}}$) forms. The result of this model inference is reported in Figures \ref{fig:step3_inferred} and \ref{fig:step3_inferredError}, illustrating the responses with respect to the raw data, the (maximal/mean) errors and eigenvalues dispersion. Note that to preserve stability one solves the structured problems \eqref{eq:ls_lODEstruct} and \eqref{eq:ls_bODEstruct} instead of  \eqref{eq:ls_lODE} and \eqref{eq:ls_bODE}. 

From Figure \ref{fig:step3_inferred}, it is clear that the inferred \textbf{ROM} models improve the accuracy of the original one. When analysing the maximal and mean errors in Figure \ref{fig:step3_inferredError}, one observes that both linear and bilinear models leads to the very good accuracy. The eigenvalue plot shows also how much the dimension number is decreased with so few eigenvalues. Moreover, both linear and bilinear realisations share the same eigenvalues in this given case. Finally, as shown on the error plots, these latter linear and bilinear \textbf{ROM} are also almost indistinguible. Indeed, this similitude was expected as the \textbf{FOM} already was able to perfectly reproduce the behaviour, with a pure linear model. Therefore, in this case, the additional bilinear term do not bring any additional accuracy improvement. This similitude between linear and bilinear is not always true, especially when a thin grid is considered for which  some nonlinear phenomena occur (see next section).

\subsection{Thin grid simulation}
\label{ssec:simulation}

\begin{figure}
    \centering
    \includegraphics[width=.8\columnwidth]{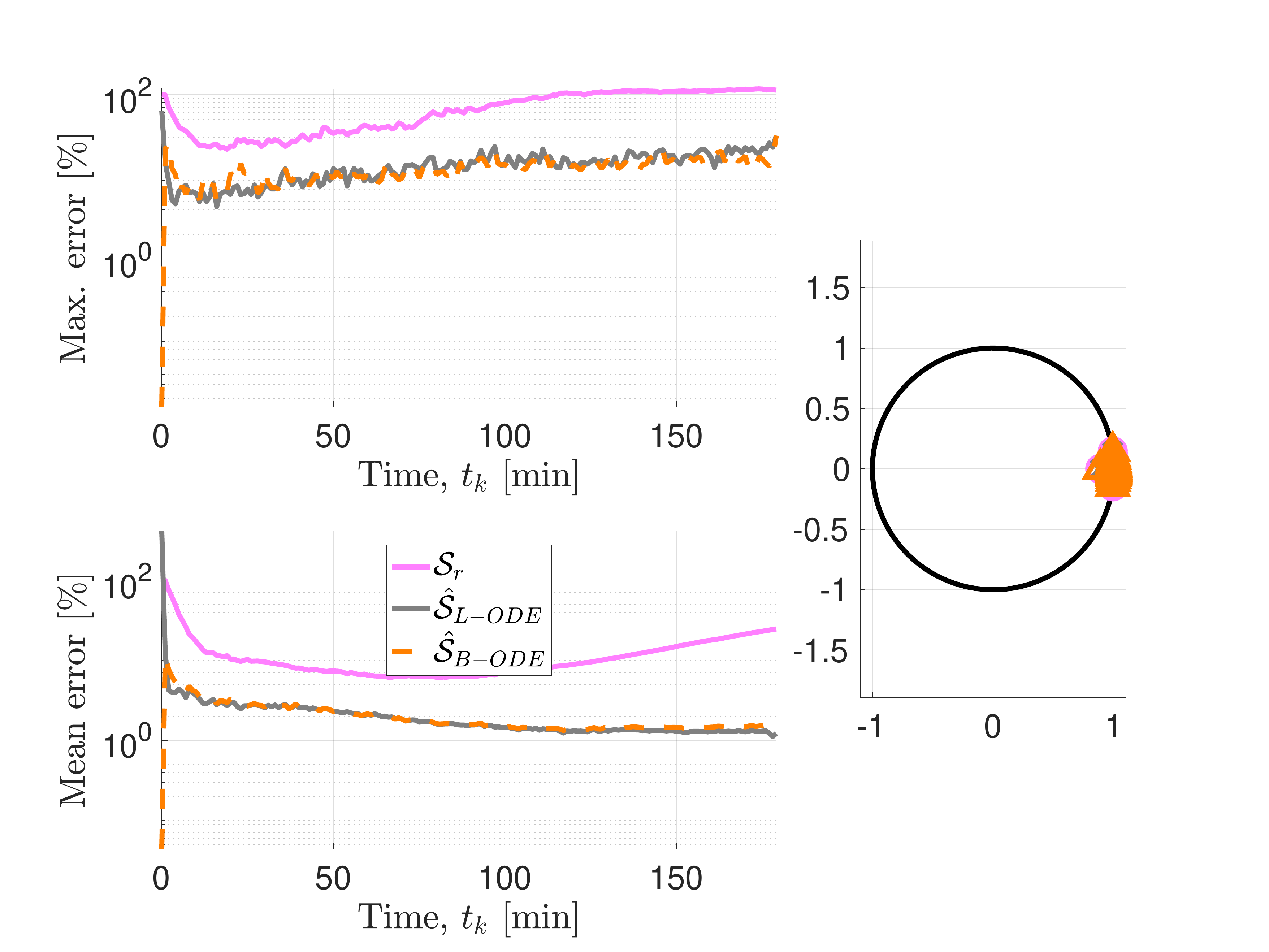}
    \includegraphics[width=.8\columnwidth]{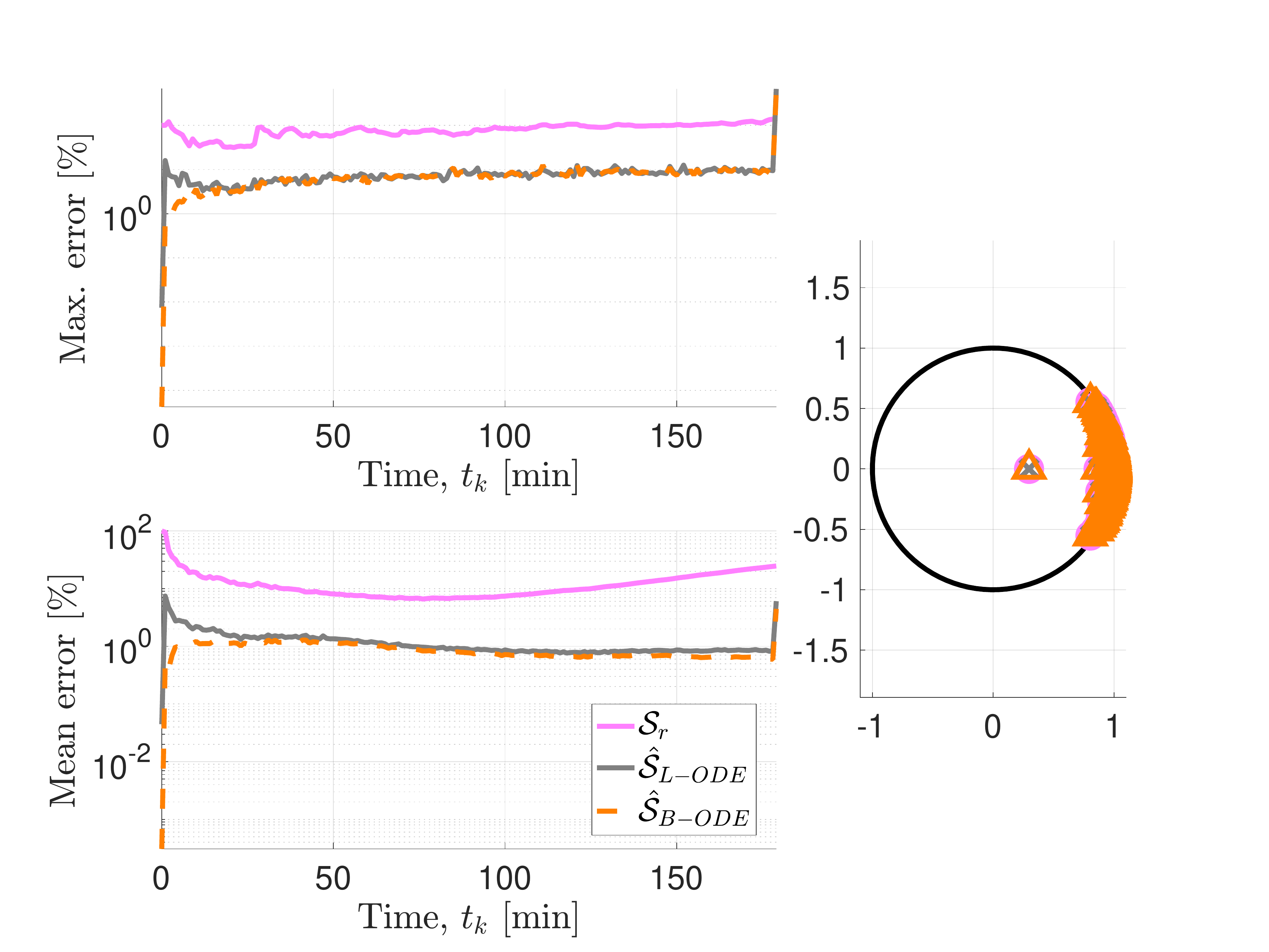}
    \caption{(Thin mesh) Top left (max) / bottom left (mean) mismatch error between the the original \textbf{LES} raw data $\mathbf Z$ and the \textbf{ROM} $\Hreal_n$ (pink), $\Hrealr_\text{L-ODE}$ (grey dotted) and $\Hrealr_\text{B-ODE}$ (orange dashed). Right: (sparse) eigenvalues of $\A_n$ (6 first ones) and $\Ar$ (Top frame $r=100$ / Bottom frame $r=150$).}
    \label{fig:inferredError_thin}
\end{figure}

Similar to the previous case, after running the \textbf{LES}, one collects the pollutants concentration $z_{i,j}(t_k)$ over each linearly spaced measurement points $\{x_i\}_{i=1}^{n_x}$ and $\{y_j\}_{j=1}^{n_y}$. Now we consider a thin $(x,y)$ spacial grid: $x_i\in[0,900]$, $y_j\in[0,900]$ (the full spacial space) with $n_x=n_y=90$ and a sampling-time $h=1$min with  $t_k=0,\dots,180$ ($N=180$).  The very same process as the above explained and shown in Figure \ref{fig:process} is now applied. Now one obtains the \textbf{FOM} $\Hreal_n$ with $n=130155$ and sparsity rate around $0.03\%$. This model is then approximated  with the stable \textbf{ROM} $\Hreal_r$, where $r=100$ and $r=250$. After model inference, one obtains the stable \textbf{ROM} $\Hrealr_\text{L-ODE}$ and $\Hrealr_\text{B-ODE}$, also of dimension $r=100$ and $r=250$. The errors are reported in Figure \ref{fig:inferredError_thin}.

\begin{figure}
    \centering
    \includegraphics[width=.45\columnwidth]{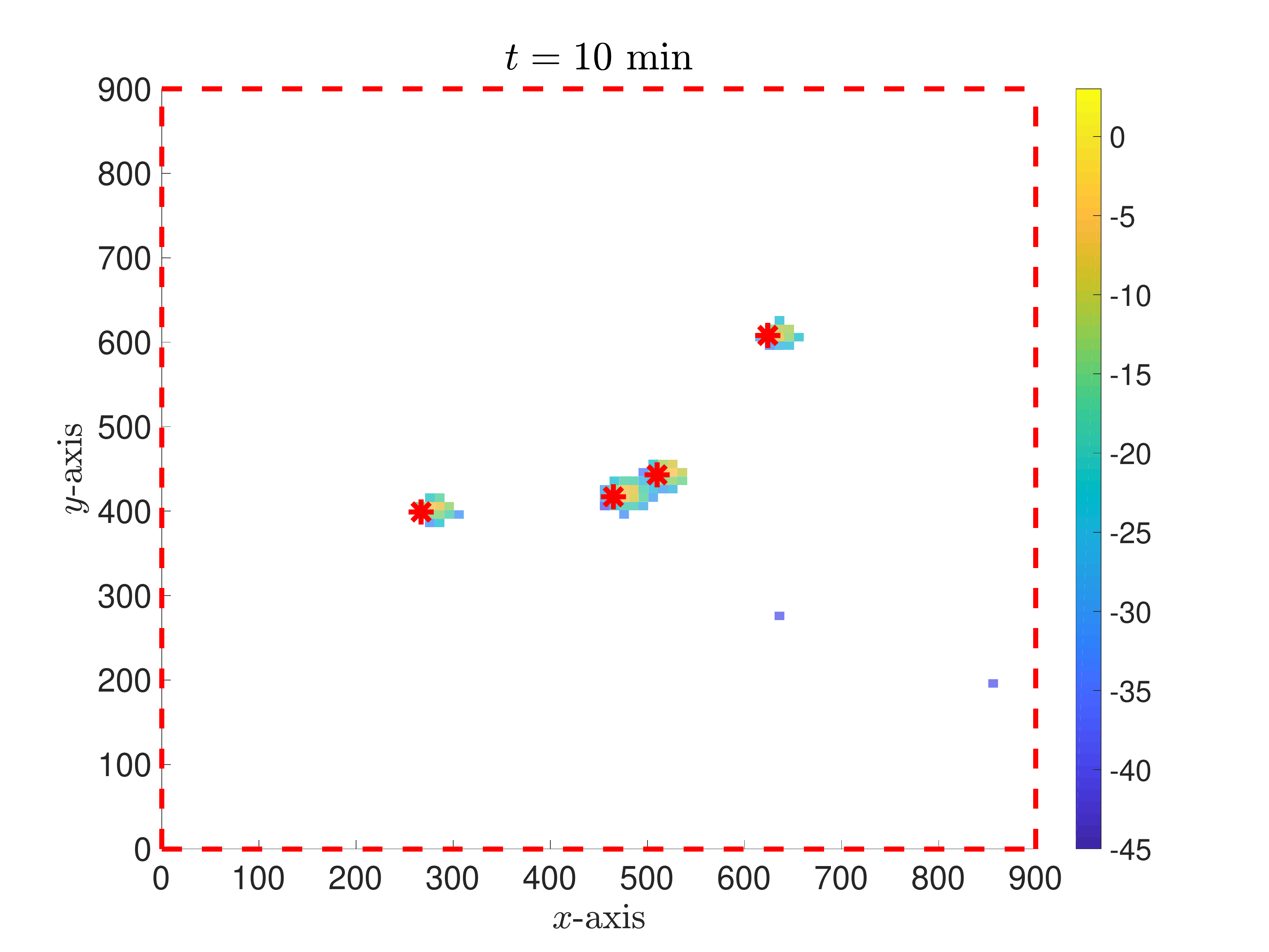}
    \includegraphics[width=.45\columnwidth]{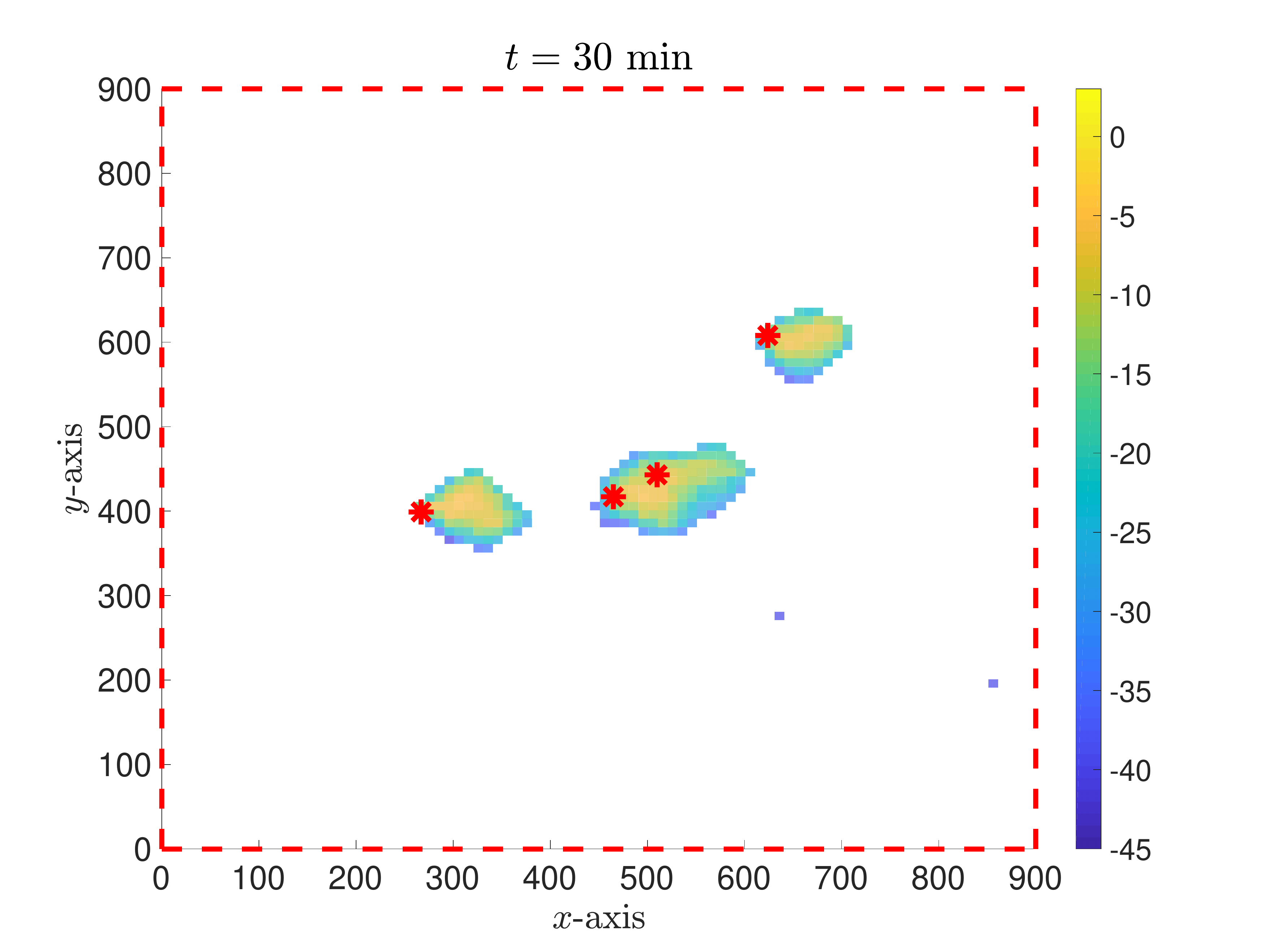}
    \includegraphics[width=.45\columnwidth]{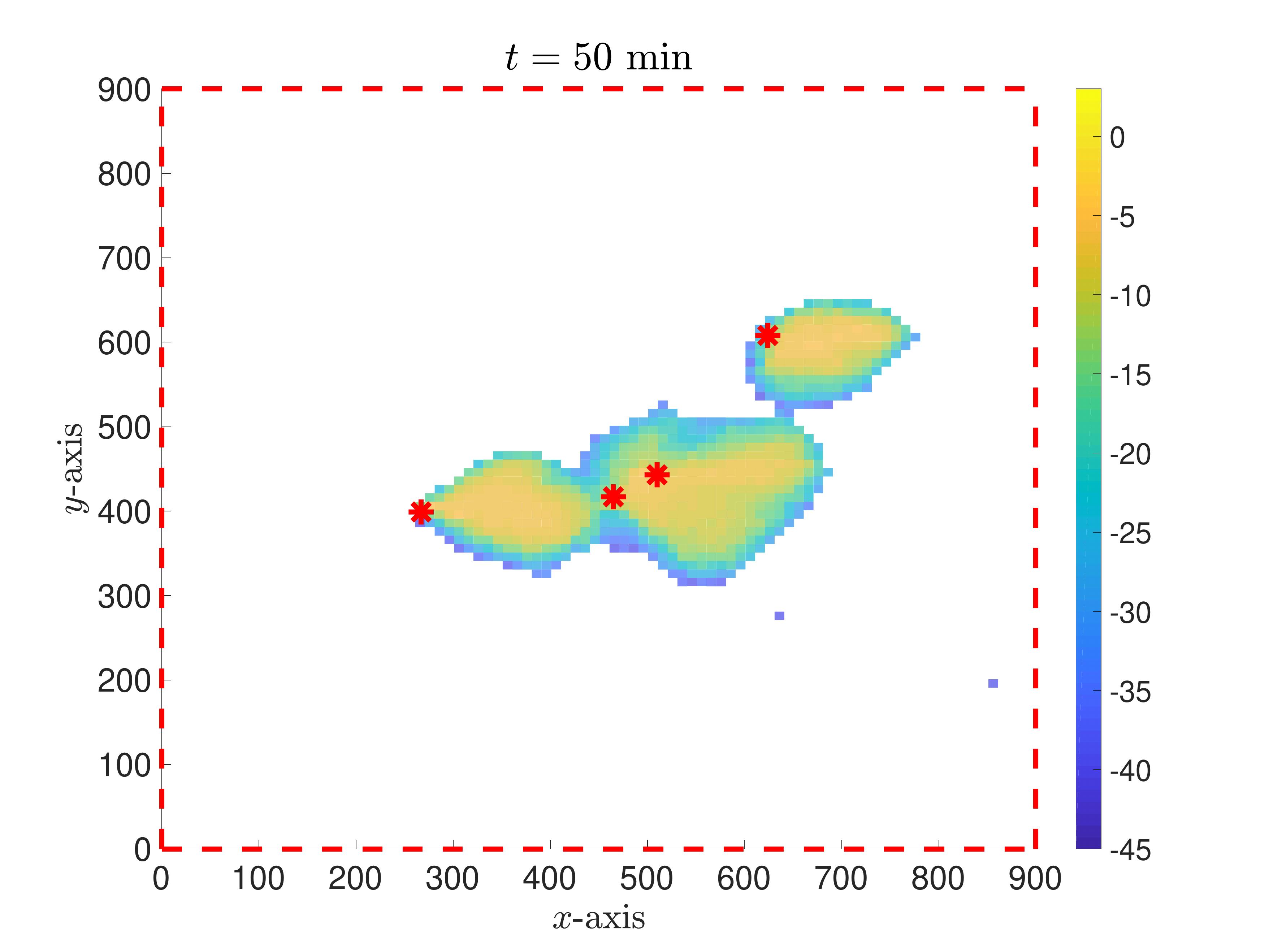}
    \includegraphics[width=.45\columnwidth]{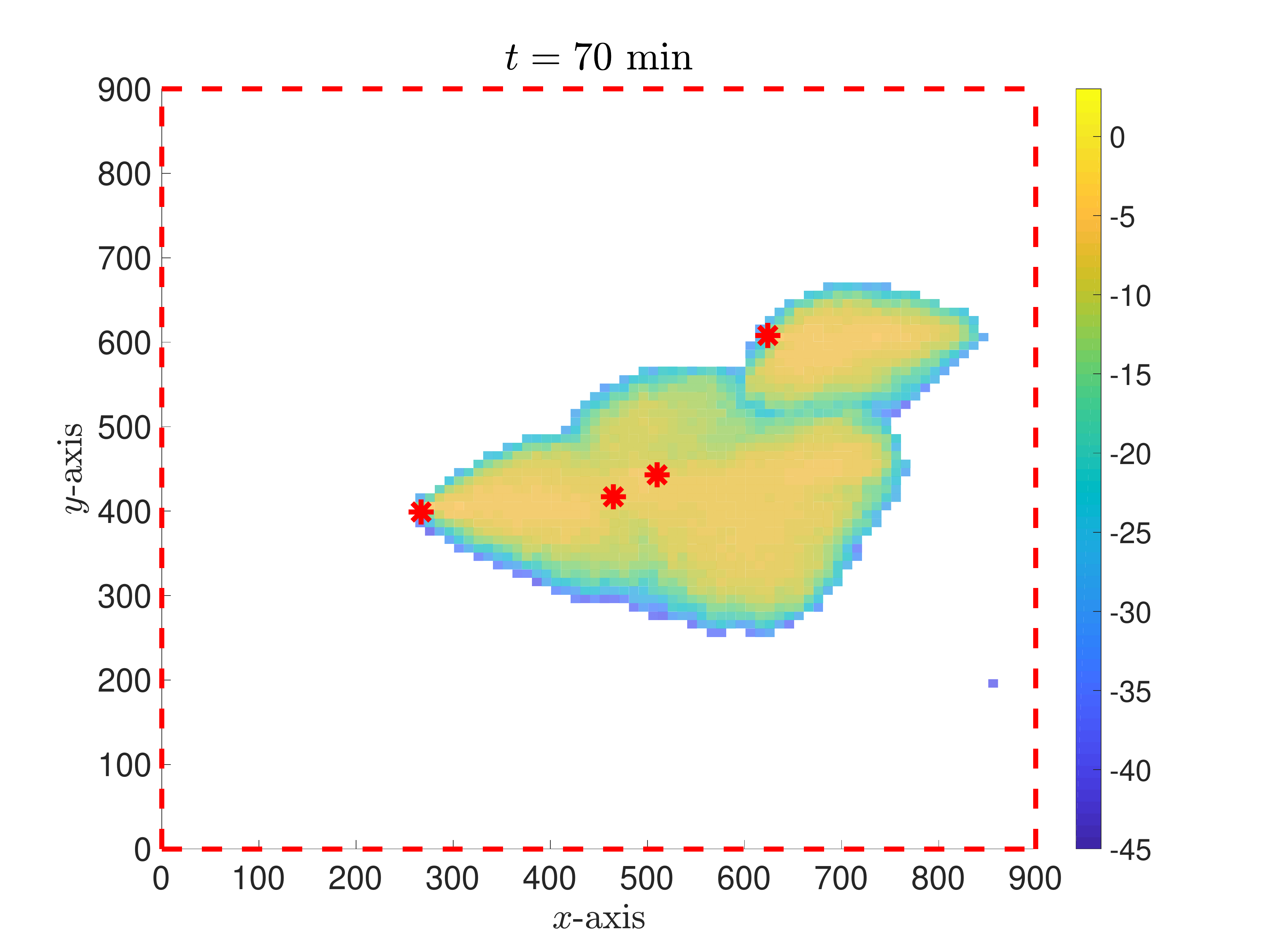}
    \includegraphics[width=.45\columnwidth]{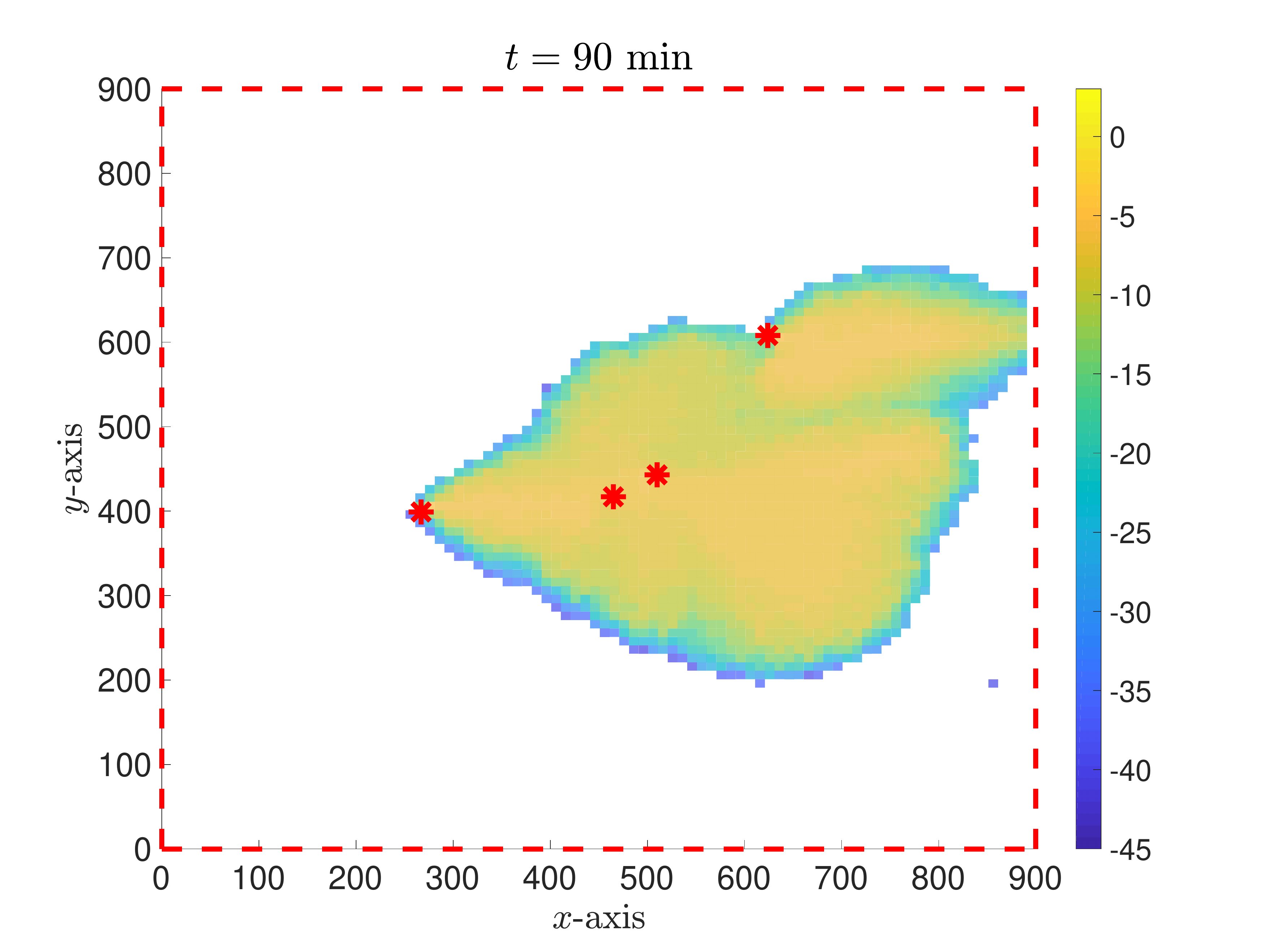}
    \includegraphics[width=.45\columnwidth]{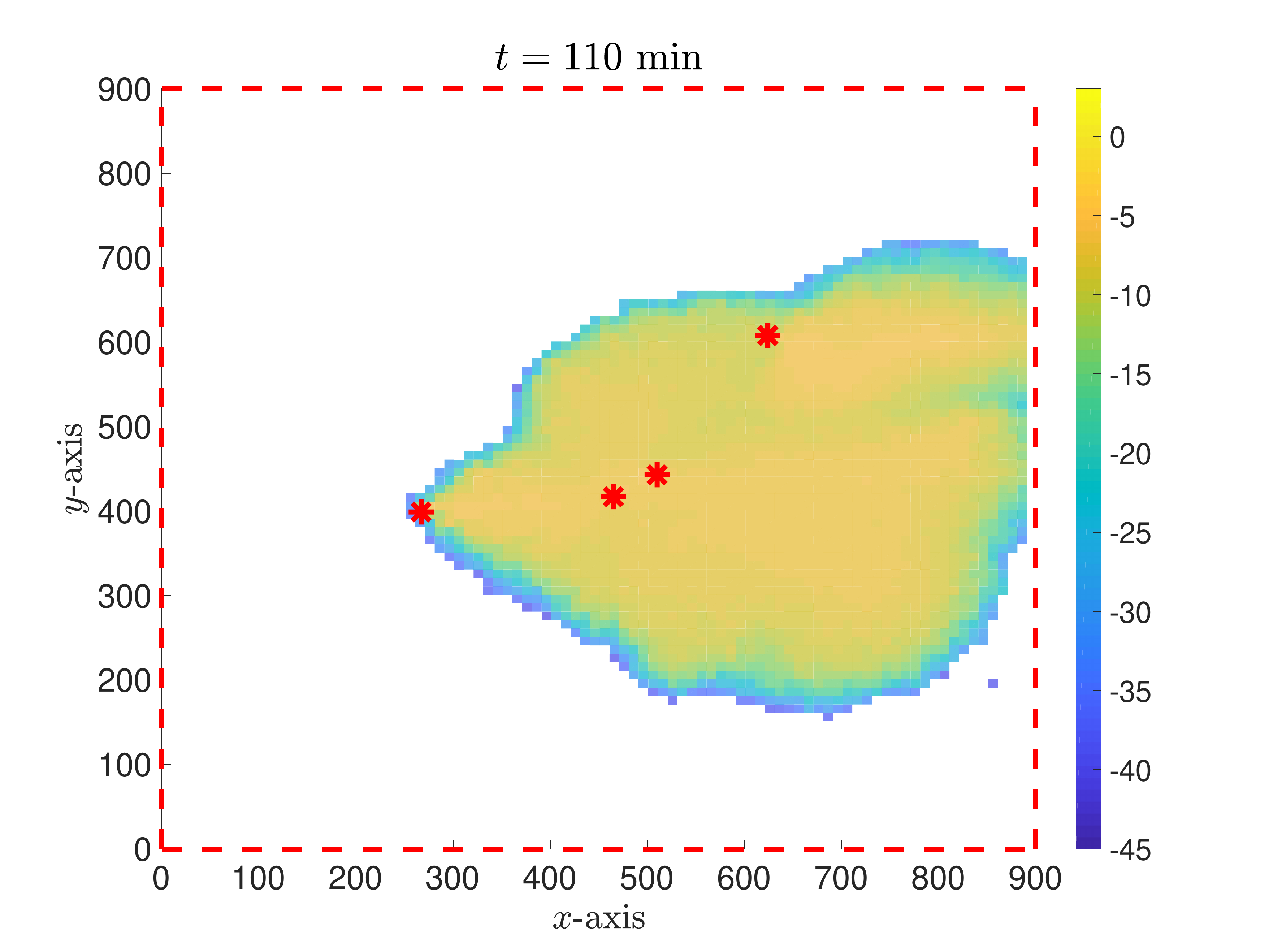}
    \caption{(Thin mesh) Snapshots of the spatial pollutants dispersion obtained with bilinear \textbf{ROM} $\Hrealr_\text{B-ODE}$ with dimension $r=250$. Pollutant values and color-bar are given with a logarithmic scale to emphasize the plume dispersion.}
    \label{fig:xyBODE}
\end{figure}

\begin{figure}
    \centering
    \includegraphics[width=.45\columnwidth]{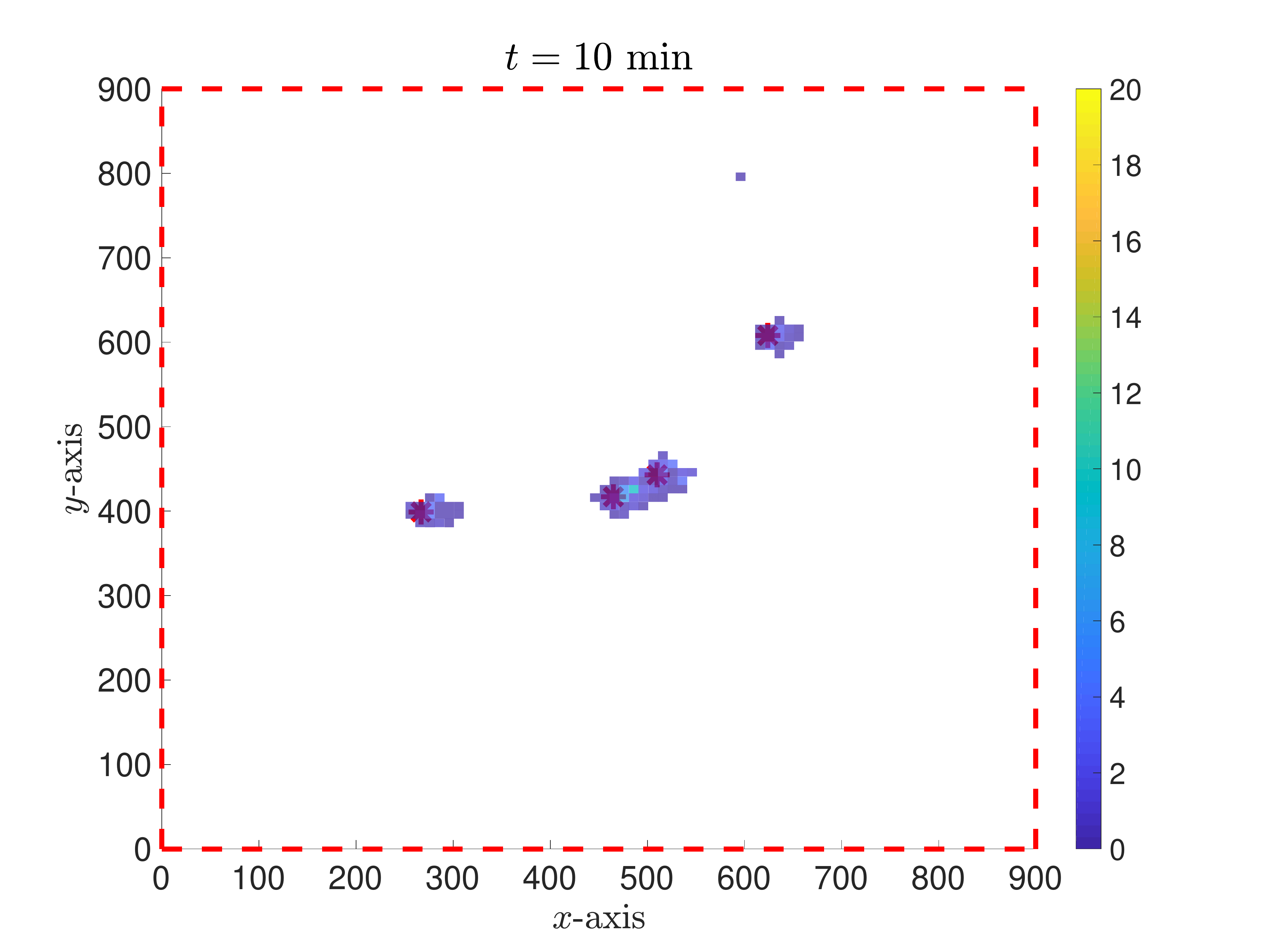}
    \includegraphics[width=.45\columnwidth]{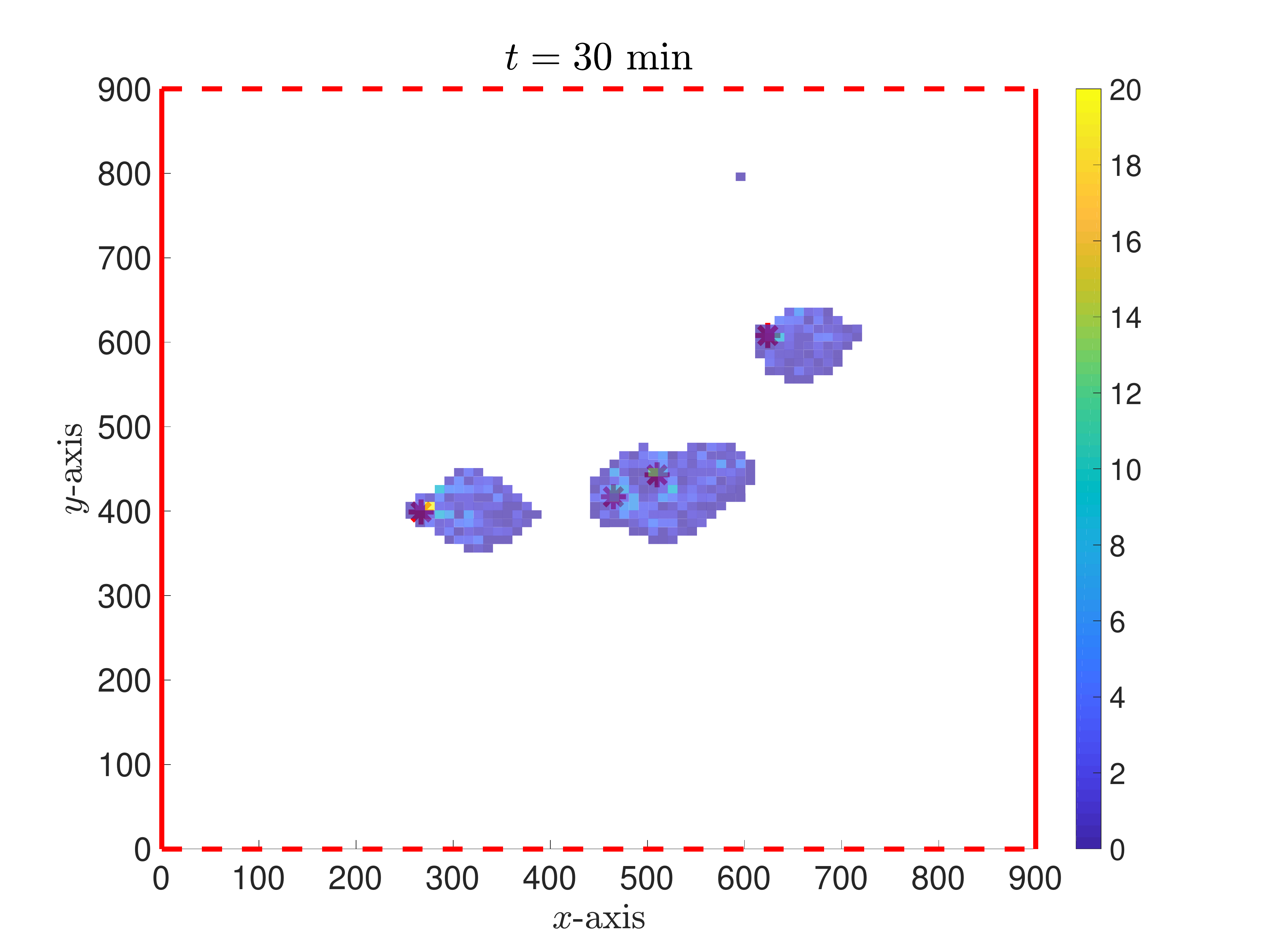}
    \includegraphics[width=.45\columnwidth]{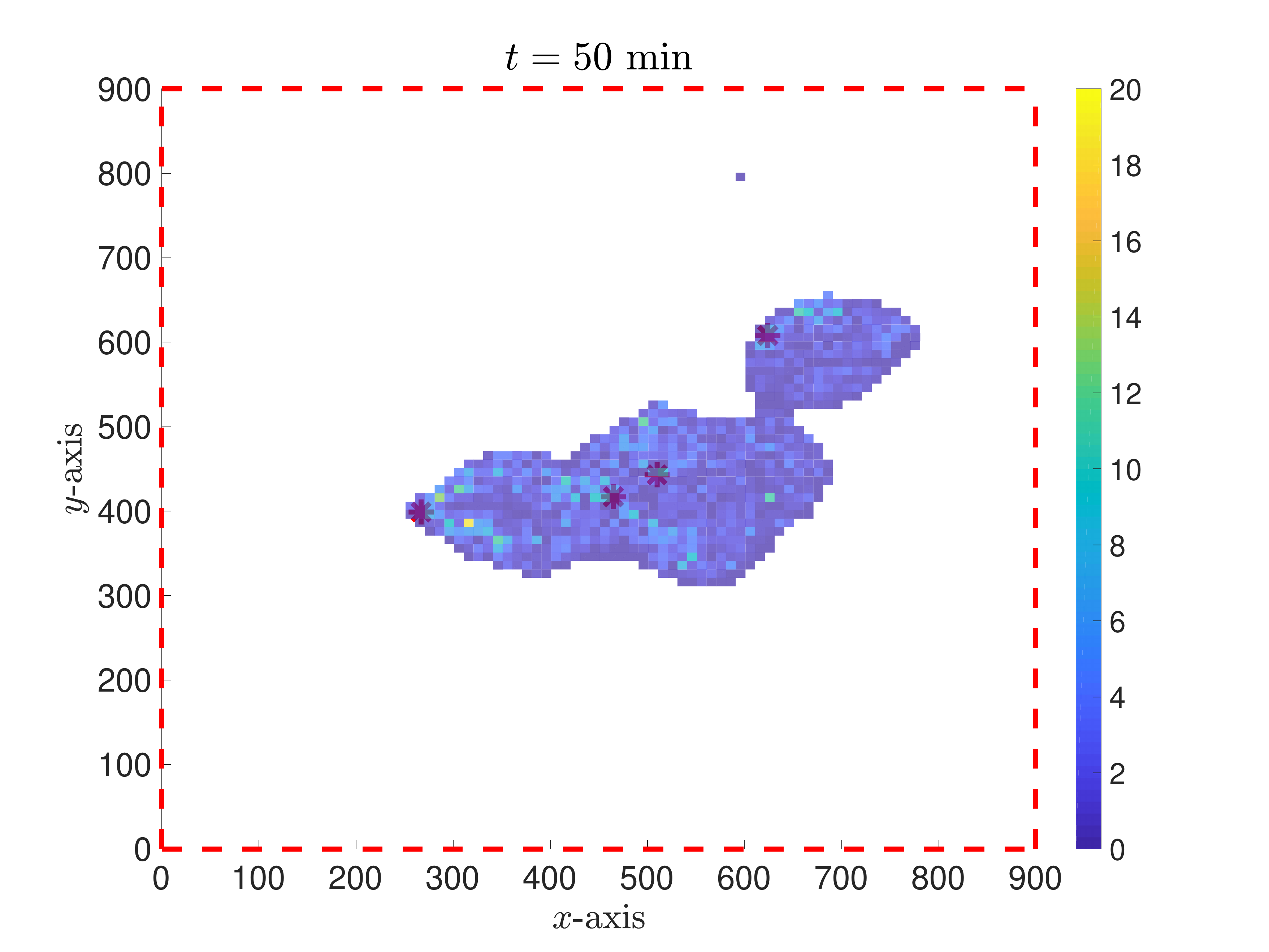}
    \includegraphics[width=.45\columnwidth]{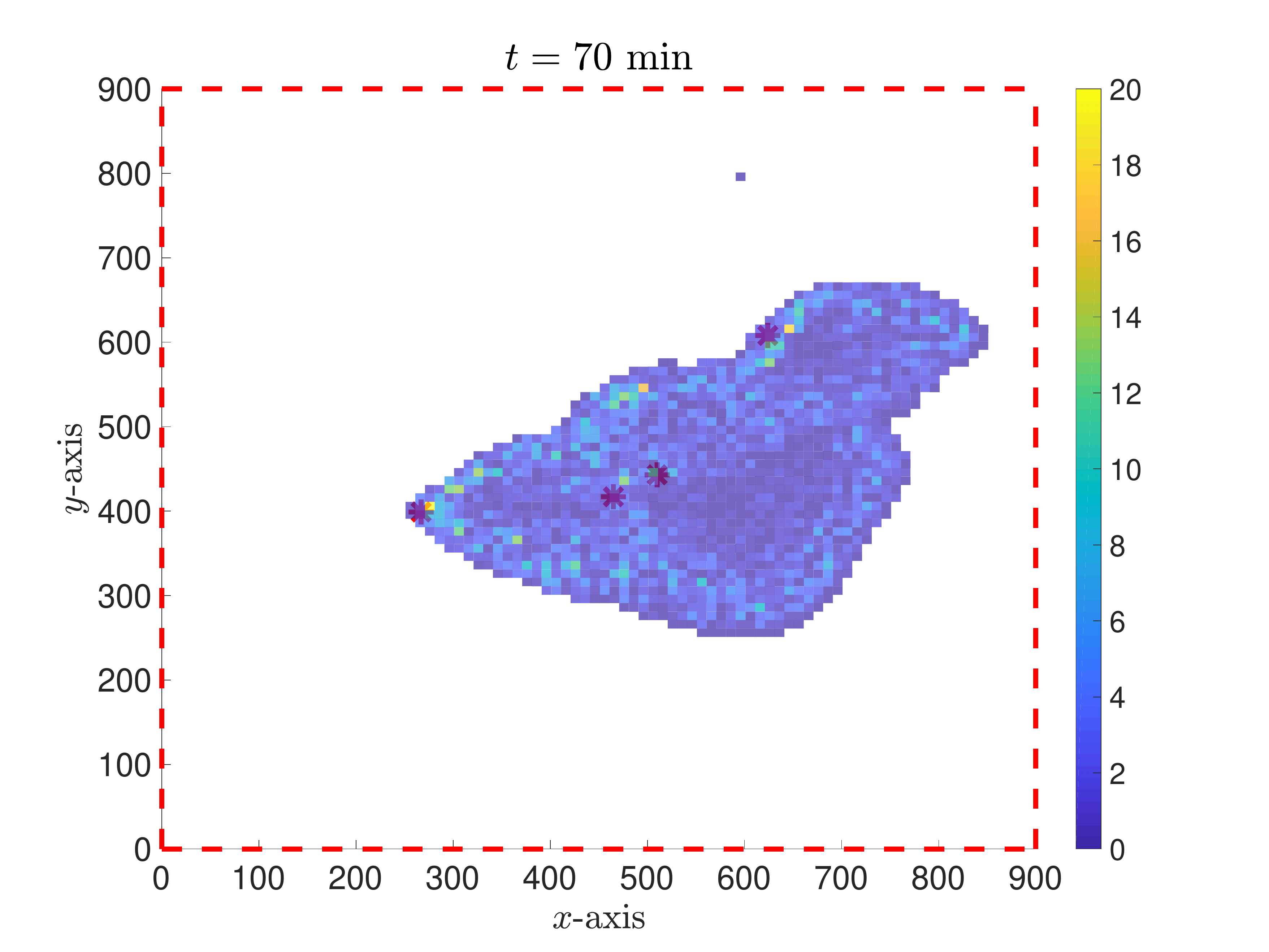}
    \includegraphics[width=.45\columnwidth]{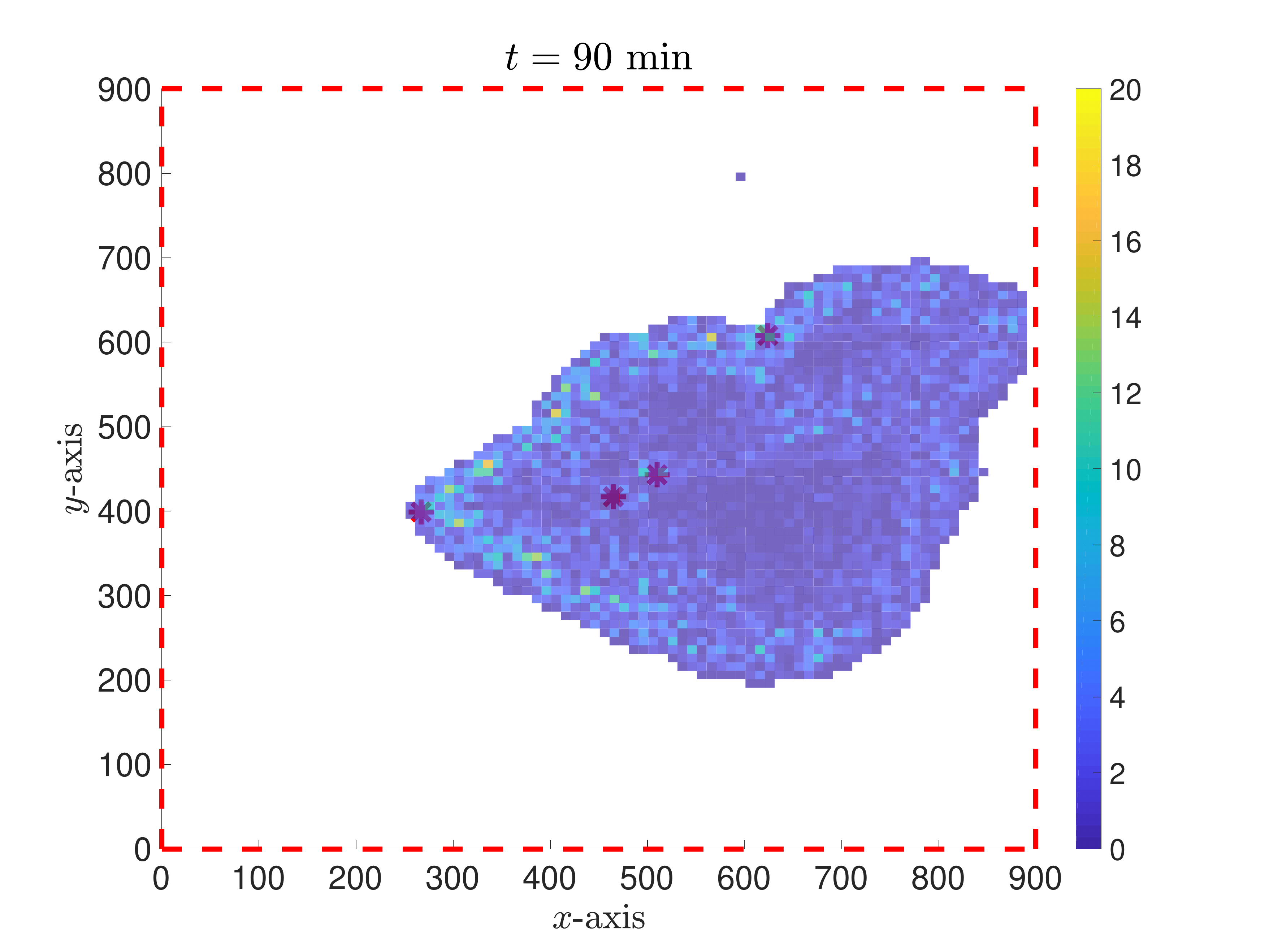}
    \includegraphics[width=.45\columnwidth]{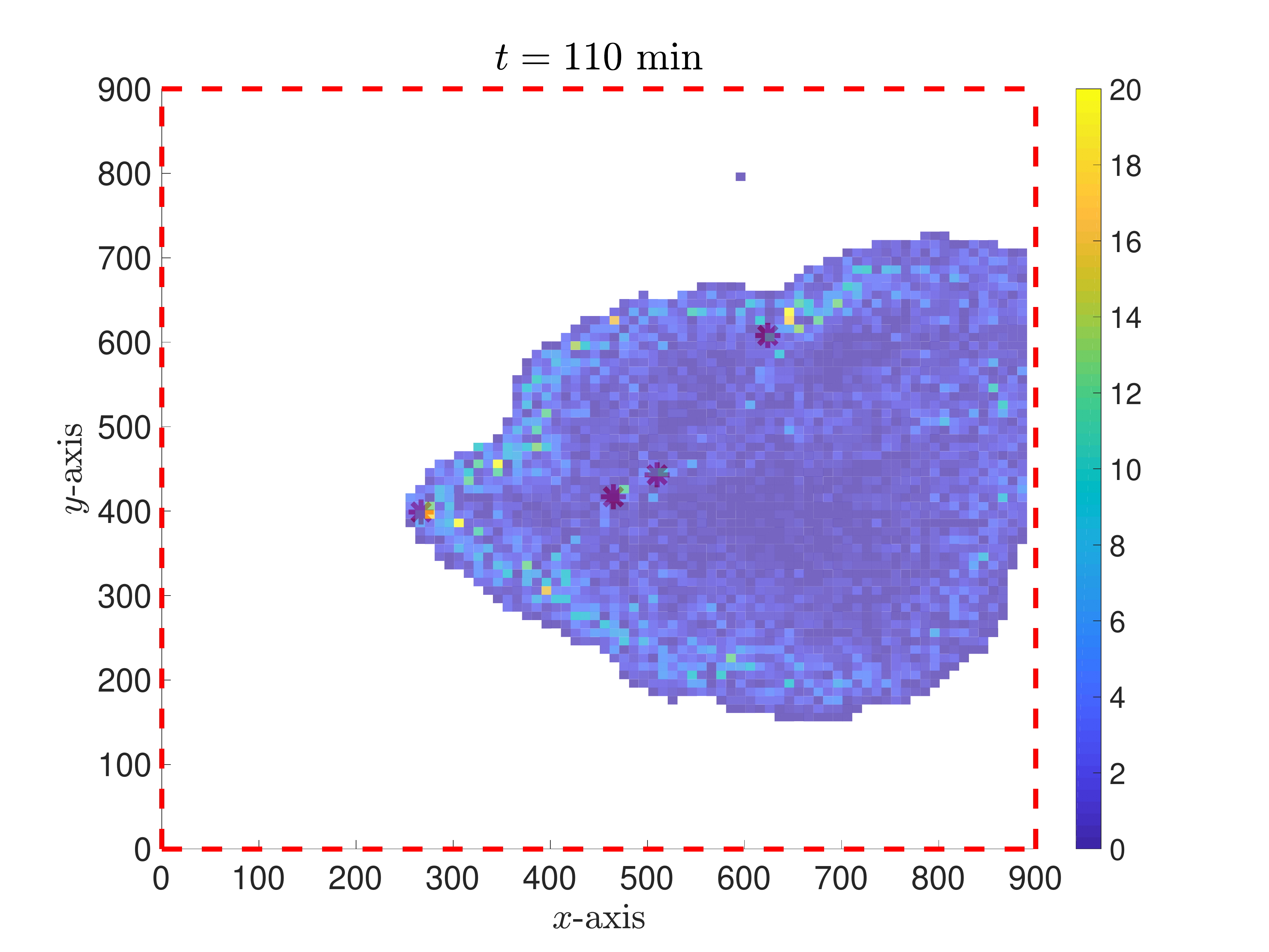}
    \caption{(Thin mesh) Snapshots of the spatial pollutants dispersion mismatch obtained with bilinear \textbf{ROM} $\Hrealr_\text{B-ODE}$ with dimension $r=250$. Mismatch error values and color-bar are given with an absolute scale as a relative percentage.}
    \label{fig:xyBODEerr}
\end{figure}

Once again, both realisations result in improving the original linear model reduction and provide a mean error around 1\% and a maximal one below 10\% which is quite reasonable considering the model simplicity. In addition, by adding the nonlinear term allows to improve the results of the linear model, especially at the beginning of the simulation where important variations occur. In this case, the nonlinear \textbf{ROM} performs better than the linear one and may be preferred. By now considering the  bilinear / state-affine model $\Hrealr_\text{B-ODE}$ model of dimension $r=250$, the complete 2-D simulation can be re-run. Horizontal cross section of pollutants concentration are illustrated through six snapshots at different time instants in Figure \ref{fig:xyBODE}. By comparing with Figure \ref{fig:xyFOM}, one appreciates the very similar restitution. In addition, the relative error computed element wise as
\begin{equation}
    \epsilon_{i,j}(t_k) = 100\frac{\abs{\hat \z_{i,j}(t_k)- \z_{i,j}(t_k)}}{\abs{\z_{i,j}(t_k)}}
\end{equation}
is displayed in percentage in Figure \ref{fig:xyBODEerr}, showing the very low residual error over the time. 

\begin{remark}[Additional material]
The proposed additional material provides a 2-D animation of the plume dispersion obtained from the original raw data collected on the complex (\textbf{LES}) simulator (left frame), the one obtained by simulating the inferred bilinear \textbf{ROM} of dimension $r=250$ (central frame), and the relative point-wise mismatch error in percentage (right frame). While the original model represents a CPU time of 5800 hours, the proposed \textbf{ROM} response is computed within a second on a standard laptop\footnote{\texttt{https://drive.google.com/file/d/1GYy8ETAFGeS9FSffbJkQ4hcBCIzrifoX/view?usp=sharing}.}.
\end{remark}

\begin{remark}[ROM computation]
In this second thin configuration, we consider a grid point of size 10. It should be noted that the \textbf{ROM} construction for this setting is around 20 minutes on a standard laptop equipped with 16Go RAM. Applying the very same approach on grid of size 1 would require additional memory to store each sub-model. This is not feasible with a standard laptop. However, the method remains valid. To bypass this limitation, a preliminary model reduction can be performed during the \textbf{FOM} construction. This is easily applicable in a practical context as it would require \eg to sub-divide the mesh and to apply local model reduction. This machinery is not exposed here but will be implemented in future developments for practical extensions.
\end{remark}

As an application, one can consider that the proposed inferred linear and nonlinear \textbf{ROM} can reasonably be used for simulation, optimisation and (Kalman) filtering for \eg prediction. Indeed, models of this complexity are affordable for multiple-query optimisation processes.


\section{Discussion and conclusions}
\label{sec:conlu}
In this paper, we first defined a completely new end to end non-intrusive time-domain reduced order nonlinear model construction method. This latter, being a mixed interpolatory - operator inference approach, stands as the principal methodological contribution of the paper. The global rationale consists in a constructive approach. First by approximating the data with a linear dynamical model, and in a second step, to adjust this latter through the addition of nonlinear terms. As it gathers the ingredients and benefits from model approximation \cite{AntoulasBook:2020} and operator inference \cite{Peherstorfer:2015}  methods, it is shown to be applicable in a very large-scale and complex context with a reasonably low computational cost (all computations are done on a standard laptop). All steps make the process easily scalable and fairly applicable to any complex simulator in a non-intrusive manner. It makes possible the construction of surrogate linear or nonlinear reduced order model from a set of input-output time-domain simulation with low human and numerical effort. 

In addition, and for the first time in meteorological simulation, the proposed approach is successfully applied on a pollutants dispersion use-case which considers very complex dynamics. The obtained surrogate model is based on input-output data, without any intrusive consideration of the multi-physics aforementioned model. 

The contributions presented in this paper open the perspective to multiple works. From the methodological and theoretical point of view related to section \ref{sec:inference}, the parametrisation of the inferred model to handle multiple configurations at the same time and to generate parametric models is an interesting path. Moreover, following the Pencil and Loewner methods philosophy, a constructive reduced order dimension and  structure procedure may represent interesting theoretical research fields to address the accuracy complexity trade-off, and to simplify user experience. From the numerical aspects, investigation of the direct \textbf{ROM} computation merging steps 2 and 3, may be an interesting topic to fasten the process and to reduce even more the computational burden and memory load.

From the application point of view, challenges cover the implementation of estimators to predict the pollutants dispersion (together with the inherent sensor placement question). Moreover, on-line model re-adjustment  from real measured data  should be considered to construct more accurate predictions to correct the inferred model by merging models and experimental data. 


\end{document}